\documentclass[11pt]{article}
\usepackage[final]{epsfig}
\usepackage{graphics}
\usepackage{amsmath}
\usepackage{amsfonts}
\usepackage{latexsym}
\usepackage{amssymb}
\usepackage{graphicx}
\usepackage{epstopdf}

\let\ssection=\section
\renewcommand{\section}{\setcounter{equation}{0}\ssection}

\newtheorem{theorem}{Theorem}
\newtheorem{lemma}{Lemma}[section]
\newtheorem{corollary}[lemma]{Corollary}

\newtheorem{proposition}[lemma]{Proposition}
\newtheorem{remark}[lemma]{Remark}
\newtheorem{example}[lemma]{Example}

\newtheorem{definition}[lemma]{Definition}

\begin{document}
\newcommand{\eps}{{\varepsilon}}
\def\a{\alpha}
\def\b{\beta}
\def\t{\tau}
\newcommand{\g}{{\gamma}}
\newcommand{\G}{{\Gamma}}
\newcommand{\proofend}{$\Box$\bigskip}
\newcommand{\C}{{\mathbb C}}
\newcommand{\Q}{{\mathbb Q}}
\newcommand{\R}{{\mathbb R}}
\newcommand{\Z}{{\mathbb Z}}
\newcommand{\N}{{\mathbb N}}
\newcommand{\RP}{{\mathbb{RP}}}
\newcommand{\cP}{{\mathcal{P}}}
\newcommand{\PSL}{{\mathrm{PSL}}}
\newcommand{\PGL}{{\mathrm{PGL}}}
\newcommand{\SL}{{\mathrm{SL}}}
\def\proof{\paragraph{Proof.}}

\newcounter{vo}
\newcommand{\vo}[1]
{\stepcounter{vo}$^{\bf VO\thevo}$%
\footnotetext{\hspace{-3.7mm}$^{\blacksquare\!\blacksquare}$
{\bf VO\thevo:~}#1}}

\newcounter{st}
\newcommand{\st}[1]
{\stepcounter{st}$^{\bf ST\thest}$%
\footnotetext{\hspace{-3.7mm}$^{\blacksquare\!\blacksquare}$
{\bf ST\thest:~}#1}}

\newcounter{rs}
\newcommand{\rs}[1]
{\stepcounter{rs}$^{\bf RS\thers}$%
\footnotetext{\hspace{-3.7mm}$^{\blacksquare\!\blacksquare}$
{\bf RS\thers:~}#1}}

\title{The Pentagram map: a discrete integrable system}

\author{Valentin Ovsienko
\and
Richard Schwartz
\and
Serge Tabachnikov}

\date{}

\maketitle

\begin{abstract}
The pentagram map is a projectively natural transformation defined
on (twisted) polygons.
A twisted polygon is a map from $\Z$ into $\RP^2$ that is
periodic modulo a projective transformation called the monodromy.
 We find a Poisson structure on the space of twisted
polygons and show that the pentagram map relative to this Poisson structure
is completely integrable.
For certain families of twisted polygons, such as those we call
universally convex, we translate the integrability
into a statement about the quasi-periodic motion for the dynamics of the
pentagram map.  We also explain how the pentagram map, in the continuous
limit, corresponds to the classical Boussinesq equation.  The Poisson
structure we attach to the pentagram map is a discrete version of the
first Poisson structure associated with the Boussinesq equation. 
A research announcement of this work appeared in \cite{OST}. 
\end{abstract}

\thispagestyle{empty}

\tableofcontents

\section{Introduction}

The notion of integrability is one of the oldest and most fundamental
notions in mathematics.
The origins of integrability lie in classical geometry and the development
of the general theory is always stimulated by the study of concrete
integrable systems.
The purpose of this paper is to study one particular dynamical
system that has a simple and natural geometric meaning and to prove
its integrability.
Our main tools are mostly geometric: the Poisson structure,
first integrals and the corresponding Lagrangian foliation.
We believe that our result opens doors for further developments
involving other approaches, such as
Lax representation, algebraic-geometric and complex analysis methods,
B\"acklund transformations;
we also expect further generalizations and relations to other fields
of modern mathematics, such as cluster algebras theory.

\begin{figure}[hbtp]
\centering
\includegraphics[height=2in]{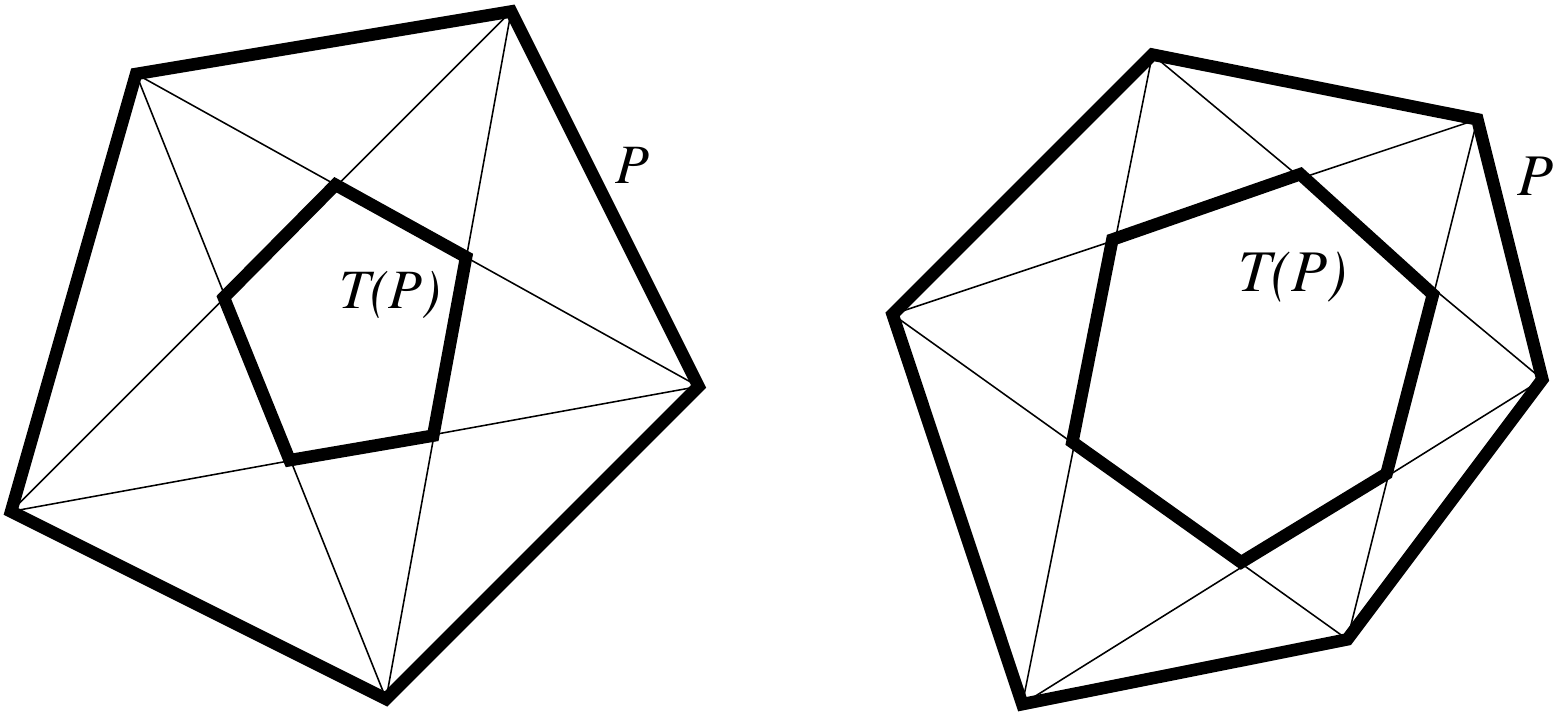}
\newline
\caption{The pentagram map defined on a pentagon and a hexagon}
\label{5and6}
\end{figure}

The {\it pentagram map}, $T$, was introduced in \cite{Sch1},
and further studied in \cite{Sch2} and \cite{Sch3}.  Originally,
the map was defined for convex closed $n$-gons.
Given such an $n$-gon $P$, the corresponding $n$-gon $T(P)$ is the convex hull of the
intersection points of consequtive shortest diagonals of $P$.
Figure~\ref{5and6} shows the situation for a convex pentagon and a convex
hexagon.  One may consider the map as defined either on unlabelled
polygons or on labelled polygons.  Later on, we shall consider the
labelled case in detail.

The pentagram  map already has some surprising features in the cases $n=5$ and $n=6$.
When $P$ is a pentagon, there is a projective transformation carrying
$P$ to $T(P)$.  This is a classical result, cf. \cite{Mot}; one of us learned
of this result from John Conway in $1987$. When $P$ is a hexagon, there is a projective
transformation carrying $P$ to $T^2(P)$.  It is not clear whether this result
was well-known to classical projective geometers, but it is easy enough
to prove.  The name pentagram map stems from the fact that
the pentagon is the simplest kind of polygon for which the map is
defined.   

Letting ${\cal C\/}_n$ denote the space of convex $n$-gons modulo projective
transformations, we can say that the pentagram map is periodic
on ${\cal C\/}_n$ for $n=5,6$.  The pentagram map certainly is not
periodic on ${\cal C\/}_n$ for $n \geq 7$.  Computer experiments
suggest that the pentagram map on ${\cal C\/}_n$ in general displays
the kind of quasi-periodic motion one sees in completely
integrable systems.  Indeed, this was conjectured (somewhat loosely) in \cite{Sch3}.
See the remarks following  Theorem 1.2 in \cite{Sch3}.

It is the purpose of this paper to establish the complete integrability
conjectured in \cite{Sch3} and to explain the underlying quasi-periodic
motion.
However, rather than work with closed $n$-gons, we will work with what we
call {\it twisted $n$-gons\/}.  A twisted $n$-gon is a map
$\phi: \Z \to \RP^2$ such that that
$$\phi(k+n)=M \circ \phi(k); \hskip 10 pt \forall k.$$
Here $M$ is some projective automorphism of $\RP^2$.  We call
$M$ the {\it monodromy\/}.  For technical reasons, we
require that every $3$ consecutive points in the image are
in general position -- i.e., not collinear.
When $M$ is the identity, we recover the notion of a closed $n$-gon.
Two twisted $n$-gons $\phi_1$ and $\phi_2$ are {\it equivalent\/}
if there is some projective transformation
$\Psi$ such that $\Psi \circ \phi_1=\phi_2$.  The two
monodromies satisfy $M_2=\Psi M_1 \Psi^{-1}$.
Let $\cP_n$ denote the space of twisted $n$-gons
modulo equivalence.

Let us emphasise that the full
space of twisted $n$-gons (rather than the geometrically natural but more restricted space of
closed $n$-gons)
is much more natural in the general context
of the integrable systems theory.
Indeed, in the ``smooth case'' it is natural to consider the full
space of linear differential equations;
the monodromy then plays an essential
r\^ole in producing  the invariants. 
This viewpoint is adopted by many authors
(see \cite{FRS, Mar} and references therein)
and this is precisely our viewpoint in the discrete case.

The pentagram map is generically defined on $\cP_n$.  However,
the lack of convexity makes it possible that the pentagram
map is not defined on some particular point of
$\cP_n$, or that the image of a point in $\cP_n$ under
the pentagram map no longer belongs to $\cP_n$.  That
is, we can lose the $3$-in-a-row property that characterizes
twisted polygons.  We will put coordinates in
$\cP_n$ so that the pentagram map becomes a rational map.
At least when $n$ is not divisible by $3$, the space
$\cP_n$ is diffeomorphic to $\R^{2n}$.  When $n$ is
divisible by $3$, the topology of the space is
trickier, but nonetheless large open subsets of
$\cP_n$ in this case are still diffeomorphic to
open subsets of $\R^{2n}$.  (Since our map is only
generically defined, the fine points of the
global topology of $\cP_n$ are not so significant.)

The action of the pentagram map in $\cP_n$ was studied
extensively in \cite{Sch3}.  In that paper, it was shown
that for every $n$ this map has a family of invariant
functions, the so-called {\it weighted monodromy invariants\/}.
There are exactly $2[n/2]+2$ algebraically independent
invariants.  Here $[n/2]$ denotes the floor of $n/2$.
When $n$ is odd, there are two exceptional monodromy functions
that are somewhat unlike the rest.  When $n$ is even,
there are $4$ such exceptional monodromy functions.
We will recall the explicit construction of these invariants in
the next section, and sketch the proofs of some of
their properties.  Later on in the paper, we shall
give a new treatment of these invariants.

Here is the main result of this paper.

\begin{theorem}
\label{main}
There exists a Poisson structure on $\cP_n$ having
co-rank $2$ when $n$ is odd and co-rank $4$ when $n$ is even.
The exceptional monodromy functions generically span
the null space of the Poisson structure, and the
remaining monodromy invariants Poisson-commute.
Finally, the Poisson structure is invariant under the pentagram map.
\end{theorem}

The exceptional monodromy functions are precisely
the {\it Casimir functions\/} for the Poisson
structure.  The generic level set of the
Casimir functions is a smooth symplectic
manifold.  Indeed, as long as we keep all
the values of the Casimir functions nonzero,
the corresponding level sets are smooth
symplectic manifolds.
The remaining monodromy invariants, when restricted to the
symplectic level sets, define a singular Lagrangian
foliation.  Generically, the dimension of the
Lagrangian leaves is precisely the same as the
number of remaining monodromy invariants.  This
is the classical picture of Arnold-Liouville
complete integrability.

As usual in this setting, the complete integrability
gives an invariant affine structure to every smooth
leaf of the Lagrangian foliation.  Relative to this
structure, the pentagram map is a translation.  Hence

\begin{corollary}
\label{closure}
Suppose that $P$ is a twisted $n$-gon that lies on a smooth
Lagrangian leaf and has a periodic orbit under the 
pentagram map.  If $P'$ is any twisted $n$-gon on the
same leaf, then $P'$ also has a periodic 
orbit with
the same period, provided that the orbit of $P'$
is well-defined.
\end{corollary}

\begin{remark}
{\rm In the result above, one can replace the word {\it periodic\/}
with $\eps$-{\it periodic\/}.  By this we mean that we
fix a Euclidean metric on the leaf and measure distances
with respect to this metric.}
\end{remark}

We shall not analyze the behavior
of the pentagram map on ${\cal C\/}_n$.
One of the difficulties in analyzing the space
${\cal C\/}_n$ of closed convex polygons modulo
projective transformations is that this space
has positive codimension in $\cP_n$ (codimension 8).  
We do not know in enough detail how the
Lagrangian singular foliation intersects
${\cal C\/}_n$, and so we cannot appeal to the
structure that exists on generic leaves.
How the monodromy invariants behave when
restricted to ${\cal C\/}_n$ is a subtle
and interesting question that we do not yet
fully know how to answer (see Theorem \ref{closed} for a partial result).
We hope to tackle the case of closed $n$-gons
in a sequel paper.

One geometric setting where our machine works perfectly is
the case of {\it universally convex $n$-gons\/}.
This is our term for a twisted $n$-gon whose image
in $\RP^2$ is strictly convex.  The monodromy
of a universally convex $n$-gon is necessarily
an element of $PGL_3(\R)$ that lifts to a
diagonalizable matrix in $SL_3(\R)$.  A universally
convex polygon essentially follows along one
branch of a hyperbola-like curve.  Let ${\cal U\/}_n$ denote
the space of universally convex $n$-gons,
modulo equivalence.  We will prove that
${\cal U\/}_n$ is an open subset of ${\cal P\/}_n$
locally diffeomorphic to $\R^{2n}$.  Further, we will see that
the pentagram map is a self-diffeomorphism of ${\cal U\/}^{2n}$. 
Finally, we will see that every
leaf in the Lagrangian foliation intersects
${\cal U\/}_n$ in a compact set.  

Combining these results with our Main Theorem and
some elementary differential topology, we arrive at the
following result.

\begin{theorem}
\label{universallyconvex}
\label{universally convex}
Almost every point of ${\cal U\/}_n$ lies on a
smooth torus that has a $T$-invariant affine
structure.  Hence, the orbit of almost every
universally convex $n$-gon undergoes quasi-periodic
motion under the pentagram map.
\end{theorem}

We will prove a variant of Theorem \ref{universally convex}
for a different family of twisted $n$-gons.  See Theorem~\ref{variant}.
The general idea is that certain points of $\cP_n$ can be interpreted
as embedded, homologically nontrivial, locally convex polygons
on projective cylinders, and suitable choices of
geometric structure give us the compactness we need
for the proof.

Here we place our results in a context. First of all, it
seems that there is some connection between our work and
cluster algebras. On the one hand,
the space of twisted polygons is known as an
example of cluster manifold, see \cite{FZ,FG} and discussion in the end of this paper.
This implies in particular that $\cP_n$ is
equipped with a canonical Poisson structure, see \cite{GSV}.
We do not know if the Poisson structure constructed in this paper
coincides with the canonical cluster Poisson structure.
On the other hand, it was shown in \cite{Sch3} that
a certain change of coordinates brings the pentagram map
rather closely in line with the octahedral recurrence, which is one of the prime examples in the theory of cluster algebras, see \cite{RR,He, Sp}.

Second of all, there is a close connection between the
pentagram map and integrable P.D.E.s.
In the last part of this paper we consider the continuous
limit of the pentagram map.
We show that this limit is precisely the classical Boussinesq equation
which is one of the best known infinite-dimensional integrable systems.
Moreover, we argue that the Poisson bracket constructed in the present
paper is a discrete analog of so-called first Poisson structure of
the Boussinesq equation.  We remark that a connection to
the Boussinesq equation was mentioned in \cite{Sch1}, but no
derivation was given.  

Discrete integrable systems is an actively developing subject, see, e.g., \cite{Ve} 
and the books \cite{Su, BS}. 
The paper \cite{BC} discusses a well-known discrete version (but with continuous time) of the Boussinesq equation; see
\cite{TN} (and references therein) for a lattice version
of this equation.
See \cite{FRS} (and references therein)
for a general theory of integrable difference equations.
Let us stress that the $r$-matrix Poisson brackets considered
in \cite{FRS} are analogous to the second
(i.e., the Gelfand-Dickey) Poisson bracket.
A geometric interpretation of all the discrete integrable systems
considered in the above references is unclear.

In the geometrical setting which is more close to our viewpoint,
see \cite{BS} for many interesting examples.
The papers \cite{Adl,Adl1} considers a discrete
integrable systems on the space of $n$-gons, different from the pentagram map.
The recent paper \cite{Kon} considers a discrete integrable
systems in the setting of projective differential geometry; some of the 
formulas in this paper are close to ours.
Finally, we mention \cite{HK,Mar} for discrete and continuous integrable systems
related both to Poisson geometry and projective differential geometry on the projective line.
\newline

We turn now to a description of the contents of the paper.
Essentially, our plan is to make a bee-line for all our
main results, quoting earlier work as much as possible.
Then, once the results are all in place, we will
consider the situation from another point of view,
proving many of the results quoted in the beginning.

One of the disadvantages of the paper \cite{Sch3}
is that many of the calculations are {\it ad hoc\/} and
done with the help of a computer.  Even though
the calculations are correct, one is not given much
insight into where they come from.  In this paper,
we derive everything in an elementary way, using
an analogy between twisted polygons and solutions
to periodic ordinary differential equations.

One might say that this paper is organized along the
lines of {\it first the facts, then the reasons\/}.
Accordingly, there is a certain redundancy in our
treatment. For instance, we introduce two natural
coordinate systems in $\cP_n$. In the first
coordinate system, which comes from \cite{Sch3},
most of the formulas are simpler.  However, the
second coordinate system, which is new, serves as
a kind of engine that drives all the derivations
in both coordinate systems; this coordinate system is 
better for computation of the monodromy too. Also, we discovered
the invariant Poisson structure by thinking
about the second coordinate system.

In \S 2 we introduce the first coordinate system,
describe the monodromy invariants, and establish
the Main Theorem.  In \S 3 we apply the main
theorem to universally convex polygons and other
families of twisted polygons.  In \S 4 we introduce
the second coordinate system.  In \S 4 and 5 we use the
second coordinate system to derive many of the
results we simply quoted in \S 2.  Finally, in \S 6
we use the second coordinate system to derive the
continuous limit of the pentagram map.

\section{Proof of the Main Theorem}

\subsection{Coordinates for the space}

In this section, we introduce our first coordinate system on the
space of twisted polygons.
As we mentioned in the introduction, a twisted $n$-gon is a
map $\phi: \Z \to \RP^2$ such that
\begin{equation}
\phi(n+k)=M \circ \phi(k)
\end{equation}
for some projective transformation $M$ and all $k$.
We let $v_i=\phi(i)$.  Thus, the vertices of
our twisted polygon are naturally $...v_{i-1},v_i,v_{i+1},...$.
Our standing assumption is that
$v_{i-1},v_i,v_{i+1}$ are in general position for all $i$,
but sometimes this assumption alone will not be
sufficient for our constructions.

The cross ratio is the most basic invariant in projective geometry.
Given four points $t_1,t_2,t_3,t_4\in\RP^1$, the cross-ratio $[t_1,t_2,t_3,t_4]$
is their unique projective invariant.
The explicit formula is as follows.
Choose an arbitrary affine parameter, then
\begin{equation}
\label{CrossR}
[t_1,t_2,t_3,t_4]=
\frac{(t_1-t_2)\,(t_3-t_4)}{(t_1-t_3)\,(t_2-t_4)}.
\end{equation}
This expression is independent of the choice of the affine parameter,
and is invariant under the 
action of $\PGL(2,\R)$ on $\RP^1$.

\begin{remark}
{\rm Many authors define the cross
ratio as the multiplicative inverse of the formula in
Equation \ref{CrossR}.  Our definition, while perhaps less
common, better suits our purposes.}
\end{remark}

The cross-ratio was used in \cite{Sch3} to define a coordinate system
on the space of twisted $n$-gons.  As the reader will see from
the definition, the construction requires somewhat more than $3$
points in a row to be in general position.  Thus, these coordinates
are not entirely defined on our space $\cP_n$.  However, they
are generically defined on our space, and this is sufficient
for all our purposes.

\begin{figure}[hbtp]
\centering
\includegraphics[height=2.7in]{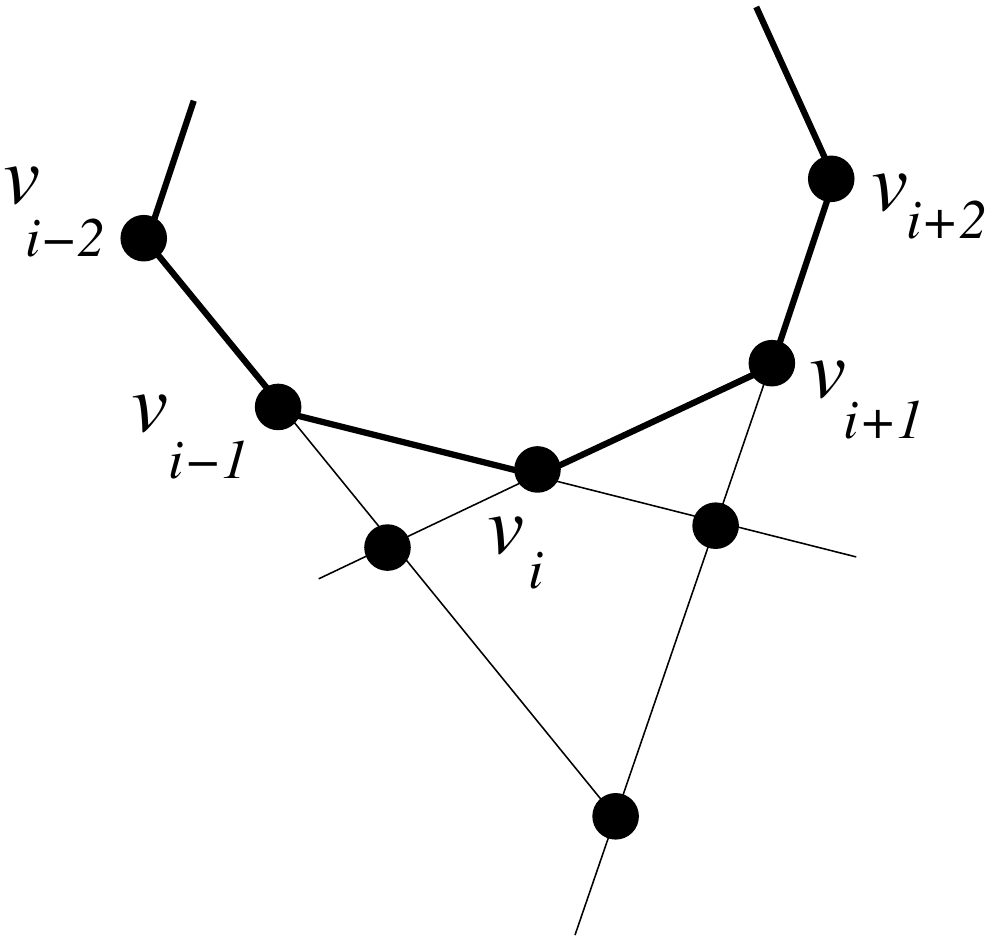}
\newline
\caption{Points involved in the definition of the invariants}
\label{corners}
\end{figure}

The construction is as follows, see Figure \ref{corners}.
We associate to every vertex $v_i$ two numbers:
\begin{equation}
\label{Corn}
\begin{array}{rcl}
x_i&=&
\displaystyle
\left[
v_{i-2},\,v_{i-1},\,
\left(
(v_{i-2},v_{i-1})\cap(v_i,v_{i+1})
\right),\,
\left(
(v_{i-2},v_{i-1})\cap(v_{i+1},v_{i+2})
\right)
\right]\\[10pt]
y_i&=&
\displaystyle
[\left(
(v_{i-2},v_{i-1})\cap(v_{i+1},v_{i+2})
\right),\,
\left(
(v_{i-1},v_{i})\cap(v_{i+1},v_{i+2})
\right),\,v_{i+1},\,v_{i+2}]
\end{array}
\end{equation}
called the left and right corner cross-ratios.
We often call our coordinates the {\it corner invariants\/}.

Clearly, the construction is $\PGL(3,\R)$-invariant and, in particular,
$x_{i+n}=x_i$ and $y_{i+n}=y_i$.
We therefore obtain a (local) coordinate system that is generically
defined on the space $\cP_n$.
In \cite{Sch3}, \S 4.2, we show how to reconstruct a twisted
$n$-gon from its sequence of invariants. 
The reconstruction is only canonical up to projective equivalence.
Thus, an attempt to reconstruct $\phi$ from $x_1,y_1,...$ 
perhaps would lead to an unequal but equivalent
twisted polygon.  This does not bother us.
The following lemma is nearly obvious.

\begin{lemma}
\label{localdiff}
At generic points, the space $\cP_n$ is locally
diffeomorphic to $\R^{2n}$.
\end{lemma}

\proof
We can perturb our sequence $x_1,y_1,...$ in any way
we like to get a new sequence
$x_1',y_1',...$.  If the perturbation is small,
we can reconstruct a new twisted $n$-gon $\phi'$
that is near $\phi$ in the following sense.  There is a
projective transformation $\Psi$ such that
$n$-consecutive vertices of $\Psi(\phi')$ are
close to the corresponding $n$ consecutive vertices
of $\phi$.  In fact, if we normalize so that
a certain quadruple of consecutive points of
$\Psi(\phi')$ match the corresponding points
of $\phi$, then the remaining points vary
smoothly and algebraically with the coordinates.
The map $(x_1',y_2',...,x_n',y_n') \to [\phi']$
(the class of $\phi'$) gives the local diffeomorphism.
\proofend

\begin{remark}
{\rm (i) Later on in the paper, we will introduce
new coordinates on all of $\cP_n$ and show, with these
new coordinates, that $\cP_n$ is globally diffeomorphic to
$\R^{2n}$ when $n$ is not divisible by $3$.  \newline
(ii)
The actual lettering we use here to define our
coordinates is different from the lettering used in
\cite{Sch3}.  Here is the correspondense:
$$...p_1,q_2,p_3,q_4... \hskip 30 pt
\Longleftrightarrow \hskip 30 pt ...,x_1,y_1,x_2,y_2,...
$$
}
\end{remark}

\subsection{A formula for the map}

\begin{figure}[hbtp]
\centering
\includegraphics[height=2.4in]{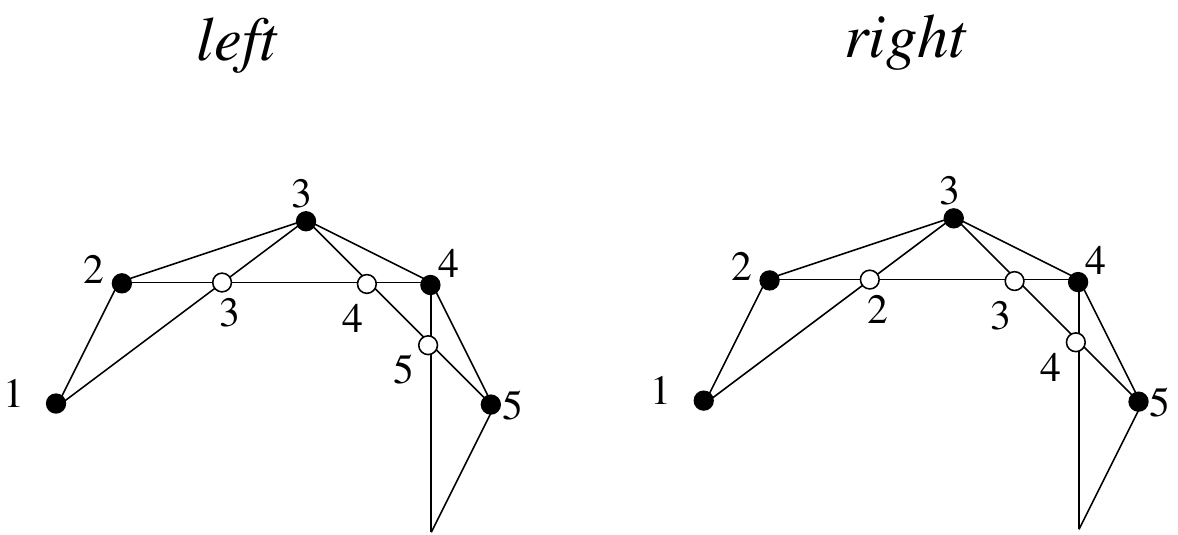}
\newline
\caption{Left and right labelling schemes}
\label{choices}
\end{figure}

In this section, we express the pentagram map in
the coordinates we have introduced in the previous section.
To save words later, we say now that we
will work with generic elements of
$\cP_n$, so that all constructions are well-defined.
Let $\phi \in \cP_n$.  Consider  the image, $T(\phi)$, of
$\phi$ under the pentagram map.  One difficulty in making
this definition is that there are two natural choices
for labelling $T(\phi)$, the {\it left choice\/} and
the {\it right choice\/}.  These choices are shown
in Figure \ref{choices}.   In the picture, the black dots represent
the vertices of $\phi$ and the white dots represent
the vertices of $T(\phi)$.  The labelling continues in
the obvious way.

If one considers the square of the pentagram map, the
difficulty in making this choice goes away.  However, 
for most of our calculations it is convenient for us
to arbitrarily choose {\it right\/} over {\it left\/}
and consider the pentagram map itself and not the square of
the map. Henceforth, we make this choice.

\begin{lemma}
\label{mapinx}
Suppose the coordinates for $\phi$ are $x_1,y_1,...$ then
the coordinates for $T(\phi)$ are
\begin{equation}
\label{ExpXEq}
T^*x_i=x_i\,\frac{1-x_{i-1}\,y_{i-1}}{1-x_{i+1}\,y_{i+1}},
\qquad 
T^*y_i=y_{i+1}\,\frac{1-x_{i+2}\,y_{i+2}}{1-x_{i}\,y_{i}},
\end{equation}
where $T^*$ is the standard pull-back of the (coordinate)
functions by the map $T$.
\end{lemma}

In \cite{Sch3}, Equation 7, we express the squared pentagram
map as the product of two involutions on $\R^{2n}$, and
give coordinates.   From this equation one can deduce
the formula in Lemma \ref{mapinx} for the pentagram map itself.
Alternatively, later in the paper we will give a self-contained
proof of Lemma \ref{mapinx}. 

Lemma \ref{mapinx} has two corollaries, which we mention here.
These corollaries are almost immediate from the formula.
First, there is an interesting scaling symmetry of the pentagram map.
We have a {\it rescaling operation\/} on $\R^{2n}$, given by the
expression
\begin{equation} 
\label{rescaling}
R_t: \hskip 10 pt (x_1,y_1,...,x_n,y_n) \to (tx_1,t^{-1}y_1,...,tx_n,t^{-1}y_n).
\end{equation}

\begin{corollary}
The pentagram map commutes with the rescaling operation.
\end{corollary}

Second, the formula for the pentagram map exhibits rather quickly
some invariants of the pentagram map.
When $n$ is odd, define
\begin{equation}
O_n=\prod_{i=1}^n x_i; \hskip 30 pt
E_n=\prod_{i=1}^n y_i.
\end{equation}
When $n$ is even, define
\begin{equation}
\label{casimir}
O_{n/2}=\prod_{i \  {\rm even\/}} x_i+
\prod_{i \  {\rm odd\/}} x_i, \hskip 20 pt
E_{n/2}=\prod_{i \  {\rm even\/}} y_i +
\prod_{i \  {\rm odd\/}} y_i.
\end{equation}
The products in this last equation run from $1$ to $n$. 

\begin{corollary}
When $n$ is odd, the functions $O_n$ and $E_n$ are
invariant under the pentagram map.
When $n$ is even, the functions $O_{n/2}$ and $E_{n/2}$ are 
also invariant under the pentagram map.
\end{corollary}

These functions are precisely the exceptional invariants we
mentioned in the introduction.  They turn out to be the
Casimirs for our Poisson structure.

\subsection{The monodromy invariants}

In this section we introduce the invariants of the pentagram
map that arise in Theorem 1.

The invariants of the pentagram map were defined and studied in \cite{Sch3}.
In this section we recall the original definition.  Later on in the
paper, we shall take a different point of view and give self-contained
derivations of everything we say here.

As above, let $\phi$ be a twisted $n$-gon with invariants $x_1,y_1,...$.  
Let $M$ be the monodromy of $\phi$.  We lift $M$ to an element
of $GL_3(\R)$.  By slightly abusing notation, we also denote this
matrix by $M$.  The two quantities
\begin{equation}
\Omega_1=\frac{{\rm trace\/}^3(M)}{{\rm det\/}(M)}; \hskip 40 pt
\Omega_2=\frac{{\rm trace\/}^3(M^{-1})}{{\rm det\/}(M^{-1})};
\end{equation}
enjoy $3$ properties.
\begin{itemize}
\item $\Omega_1$ and $\Omega_2$ are independent of the lift of $M$.
\item $\Omega_1$ and $\Omega_2$ only depend on the conjugacy class of $M$.
\item $\Omega_1$ and $\Omega_2$ are rational functions in the corner invariants.
\end{itemize}

We define
\begin{equation}
\widetilde \Omega_1=O_n^2E_n \Omega_1; \hskip 30 pt
\widetilde \Omega_2=O_nE_n^2 \Omega_2.
\end{equation}
In \cite{Sch3} it is shown that
$\widetilde \Omega_1$ and $\widetilde \Omega_2$ are 
polynomials in the corner invariants.
Since the pentagram map preserves the monodromy, and 
$O_n$ and $E_n$ are invariants, the two functions
$\widetilde \Omega_1$ and $\widetilde \Omega_2$
are also invariants.

We say that a polynomial in the corner invariants
has {\it weight\/} $k$ if we have the following
equation
\begin{equation}
R_t^*(P)=t^k P.
\end{equation}
here $R_t^*$ denotes the natural operation on polynomials
defined by the rescaling operation (\ref{rescaling}).  For instance,
$O_n$ has weight $n$ and $E_n$ has weight $-n$.  In
\cite{Sch3} it is shown that
\begin{equation} \label{OandE}
\widetilde \Omega_1=\sum_{k=1}^{[n/2]} O_k; \hskip 30 pt
\widetilde \Omega_2=\sum_{k=1}^{[n/2]} E_k
\end{equation}
where $O_k$ has weight $k$ and $E_k$ has weight $-k$.
Since the pentagram map commutes with the rescaling
operation and preserves $\widetilde \Omega_1$ and
$\widetilde \Omega_2$, it also preserves their
``weighted homogeneous parts''.   That is, the
functions $O_1,E_1,O_2,E_2,...$ are also invariants
of the pentagram map.  These are the monodromy
invariants.  They are all nontrivial polynomials.
\newline
\newline
{\bf Algebraic Independence:\/}
In \cite{Sch3}, \S 6, it is shown
that the monodromy
invariants are algebraically independent
provided that, in the even case, we ignore
$O_{n/2}$ and $E_{n/2}$. 
We will not reproduce the proof in this paper,
so here we include a brief description of the
argument. Since we are mainly trying to give the
reader a feel for the argument, we will explain
a variant of the method in \cite{Sch3}.
Let $f_1,...,f_k$ be the complete list of
invariants we have described above.  Here
$k=2[n/2]+2$.  
  If our functions
were not algebraically independent, then 
the gradients $\nabla f_1,...,\nabla f_k$
would never be linearly independent.
To rule this out, we just have to establish
the linear independence at a single point.
One can check this at the point
$(1,\omega,...,\omega^{2n})$, where $\omega$ is a
$(4n)$th root of unity.   The actual method
in \cite{Sch3} is similar to this, but uses
a trick to make the calculation easier.
Given the formulas for the invariants we present
below, this calculation is really just a matter
of combinatorics.   Perhaps an easier calculation
can be made for the point $(0,1,...,1)$, which
also seems to work for all $n$.

\subsection{Formulas for the invariants}

In this section, we recall the explicit formulas for the monodromy
invariants given in \cite{Sch3}.
Later on in the paper, we will give a self-contained
derivation of the formulas.  From the point
of view of our main theorems, we do not need to know
the formulas, but only their algebraic independence
and Lemma \ref{FirstClaim} below.

We introduce the monomials
\begin{equation}
\label{uniT}
X_i:=x_i\,y_i\,x_{i+1}.
\end{equation}
\begin{enumerate}
\item
We call two monomials $X_i$ and $X_j$ consecutive if 
$
j\in\left\{i-2,\,i-1,\,i,\,i+1,\,i+2\right\};
$
\item
we call $X_i$ and $x_j$ consecutive if
$j\in\left\{i-1,\,i,\,i+1,\,i+2\right\};$
\item
we call $x_i$ and $x_{i+1}$ consecutive.
\end{enumerate}

Let $O(X,x)$ be a monomial obtained by the product of the monomials
$X_i$ and $x_j$, i.e.,
$$
O=X_{i_1}\cdots{}X_{i_s}\,x_{j_1}\cdots{}x_{j_t}.
$$
Such a monomial is called admissible if no two of the indices 
are consecutive.
For every admissible monomial, we define the weight $|O|$
and the sign $\mathrm{sign}(O)$ by
$$
|O|:=s+t,
\qquad
\mathrm{sign}(O):=(-1)^t.
$$

With these definitions, it turns out that
\begin{equation}
\label{InvOk}
O_k=\sum_{|O|=k}\mathrm{sign}(O)\,O; \hskip 30 pt k\in\left\{1,2,\ldots,\left[\frac{n}{2}\right]\right\}
\end{equation}
The same formula works for $E_k$, if we make all the same
definitions with $x$ and $y$ interchanged.

\begin{example}
{\rm
For $n=5$ one obtains the following polynomials
$$
O_1=\sum_{i=1}^5
\left(
x_i\,y_i\,x_{i+1}-x_i
\right),
\qquad
O_2=\sum_{i=1}^5
\left(
x_i\,x_{i+2}-
x_i\,y_i\,x_{i+1}\,x_{i+3}
\right)
$$
together with $O_5$.
}
\end{example}

Now we mention the needed symmetry property.
Let $\tau$ be the involution on the indices: 
\begin{equation}
\label{Sig}
\tau:{}
\left\{
\begin{array}{l}
x_i \mapsto x_{1-i}\\
y_i \mapsto y_{-i}
\end{array}
\right.
\qquad
\hbox{mod}\, n
\end{equation}
Then $\tau$ acts on the variables, monomials and polynomials.
\smallskip

\begin{lemma}
\label{FirstClaim}
One has $\tau(O_k) = O_k$.
\end{lemma}

\proof
$\tau$ takes an admissible partition to an admissible one and does not change the number of singletons involved.
\proofend

\subsection{The Poisson bracket}

In this section, we introduce the Poisson bracket on $\cP_n$.
Let $C^{\infty}_n$ denote the algebra of
smooth functions on $\R^{2n}$.  A
{\it Poisson structure\/} on $C^{\infty}_n$ is a
map
\begin{equation}
\{\ ,\ \}: C^{\infty}_n \times C^{\infty}_n \to C^{\infty}_n
\end{equation}
that obeys the following axioms.
\begin{enumerate}
\item Antisymmetry: $\{f,g\}=-\{g,f\}$
\item Linearity: $\{af_1+f_2,g\}=a\{f_1,g\}+\{f_2,g\}$.
\item Leibniz Identity: $\{f,g_1g_2\}=g_1\{f,g_2\}+g_2\{f,g_1\}$.
\item Jacobi Identity: $\Sigma \{f_1,\{f_2,f_3\}\}=0$.
\end{enumerate}
Here $\Sigma$ denotes the cyclic sum.

We define the following Poisson bracket on the coordinate
functions of $\R^{2n}$:
\begin{equation} \label{PB}
\{x_i, x_{i\pm1}\}=\mp{}x_i\,x_{i+1},
\qquad
\{y_i,y_{i\pm1}\}=\pm{}y_i\,y_{i+1}.
\end{equation}
All other brackets not explicitly mentioned above vanish.
For instance $$\{x_i,y_j\}=0; \hskip 30 pt \forall\ i,j.$$
Once we have the definition on the coordinate functions,
we use linearity and the Liebniz rule to extend to all
rational functions.  Though it is not necessary for
our purposes, we can extend to all smooth functions by
approximation.  Our formula automatically builds in the
anti-symmetry. Finally,
for a ``homogeneous bracket'' as we have defined, it is well-known
(and an easy exercise) to show that the Jacobi identity holds.

Henceforth we refer to {\it the Poisson bracket\/} as the one
that we have defined above.
Now we come to one of the central results in the paper.
This result is our main tool for establishing the
complete integrability and the quasi-periodic motion.

\begin{lemma}
\label{invbrack}
The Poisson bracket is invariant with respect to the pentagram map.
\end{lemma}

\proof
Let $T^*$ denote the action of the pentagram map on rational
functions. 
One has to prove that for any two functions $f$ and $g$ one has
$\{T^*(f),T^*(g)\}=\{f,g\}$ and of course 
it  suffices to check this fact for the coordinate functions.
We will use the explicit formula (\ref{ExpXEq}).

To simplify the formulas, we introduce the following notation:
$\varphi_i=1-x_i\,y_i$.
Lemma \ref{mapinx} then reads:
$$
T^*(x_i)=x_i\,\frac{\varphi_{i-1}}{\varphi_{i+1}},
\qquad
T^*(y_i)=y_{i+1}\,\frac{\varphi_{i+2}}{\varphi_{i}}.
$$
One easily checks that
$
\{\varphi_i,\varphi_j\}=0
$
for all $i,j$. 
Next,
$$
\begin{array}{rcl}
\{x_i,\varphi_j\}&=&
\displaystyle
\left(
\delta_{i,j-1}-\delta_{i,j+1}
\right)
x_i\,x_j\,y_j\\[5pt] 
\{y_i,\varphi_j\}&=&
\displaystyle
\left(
\delta_{i,j+1}-\delta_{i,j-1}
\right)
x_j\,y_i\,y_j.
\end{array} 
$$

In order to check the $T$-invariance of the bracket,
one has to check that the relations between the functions
$T^*(x_i)$ and $T^*(y_j)$ are the same as for $x_i$ and $y_j$.
The first relation to check is: $\{T^*(x_i),T^*(y_j)\}=0$.

Indeed,
$$
\begin{array}{rcl}
\left\{T^*(x_i),T^*(y_j)\right\}&=&
\displaystyle
\{x_i,\varphi_{j+2}\}\,
\frac{y_{j+1}\,\varphi_{i-1}}{\varphi_{i+1}\,\varphi_j} - 
\{x_i,\varphi_j\}\,
\frac{y_{j+1}\,\varphi_{i-1}\,\varphi_{j+2}}{\varphi_{i+1}\,\varphi_j^2}\\[16pt]
&&
\displaystyle
-\{y_{j+1},\varphi_{i-1}\}\,
\frac{x_i\,\varphi_{j+2}}{\varphi_{i+1}\,\varphi_j}+
\{y_{j+1},\varphi_{i+1}\}\,
\frac{x_i\,\varphi_{i-1}\,\varphi_{j+2}}{\varphi_{i+1}^2\,\varphi_j}\\[16pt]
&=&
\displaystyle
\left(
\delta_{i,j+1}-\delta_{i,j+3}
\right)
\frac{x_i\,x_{j+2}\,y_{j+2}\,y_{j+1}\,\varphi_{i-1}}{\varphi_{i+1}\,\varphi_j}\\[16pt]
&&
\displaystyle
-\left(
\delta_{i,j-1}-\delta_{i,j+1}
\right)
\frac{x_i\,x_j\,y_j\,y_{j+1}\,\varphi_{i-1}\,\varphi_{j+2}}{\varphi_{i+1}\,\varphi_j^2}\\[16pt]
&&
\displaystyle
-\left(
\delta_{j+1,i}-\delta_{j+1,i-2}
\right)
\frac{x_{i-1}\,y_{j+1}\,y_{i-1}\,x_i\,\varphi_{j+2}}{\varphi_{i+1}\,\varphi_j}\\[16pt]
&&
\displaystyle
+\left(
\delta_{j+1,i+2}-\delta_{j+1,i}
\right)
\frac{x_{i+1}\,y_{j+1}\,y_{i+1}\,x_i\,\varphi_{i-1}\,\varphi_{j+2}}{\varphi_{i+1}^2\,\varphi_j}\\[16pt]
&=&
0,
\end{array}
$$
since the first term cancels with the third and the second with the last one.

One then computes $\{T^*(x_i),T^*(x_j)\}$ and $\{T^*(y_i),T^*(y_j)\}$,
the computations are similar to the above one and will be omitted.
\proofend

Two functions $f$ and $g$ are said to {\it Poisson commute\/} if $\{f,g\}=0$.

\begin{lemma}
The monodromy invariants Poisson commute.
\end{lemma}

\proof 
Let $\tau$ by the involution on the indices
defined at the end of the last section.
We have $\tau(O_k) = O_k$ by Lemma \ref{FirstClaim}.
We make the following claim:
For all polynomials $f(x,y)$ and $g(x,y)$, one has 
$$
\{\tau(f),\tau(g)\}=-\{f,g\}.
$$
Assuming this claim, we have
$$
\{O_k, O_l\}=\{\tau(O_k), \tau(O_l)\}=-\{O_k, O_l\},
$$
hence the bracket is zero.
The same argument works for $\{E_k,E_l\}$ and $\{O_k,E_l\}$.

Now we prove our claim.
It suffices to check the claim when  $f$ and $g$ are monomials in variables $(x,y)$.
In this case, we have: $\{f,g\}=C fg$ where $C$ is the sum of $\pm 1$, 
corresponding to ``interactions" between factors $x_i$ in $f$ and $x_j$
in $g$ (resp. $y_i$ and $y_j$).
Whenever a factor $x_i$ in $f$ interacts with a factor $x_j$ in $g$
(say, when $j=i+1$, and the contribution is +1), 
there will be an interaction of $x_{-i}$ in $\tau(f)$ and $x_{-j}$ in 
$\tau(g)$ yielding the opposite sign 
(in our example, $-j=-i-1$, and the contribution is $-1$).  This
establishes the claim, and
hence the lemma.
\proofend

A function $f$ is called a {\it Casimir\/} for the Poisson bracket if
$f$ Poisson commutes with all other functions.  It suffices to check
this condition on the coordinate functions.  An easy calculation
yields the following lemma.  We omit the details.

\begin{lemma}
\label{cas}
The invariants in Equation \ref{casimir} are Casimir functions
for the Poisson bracket.
\end{lemma}

\subsection{The corank of the structure}

In this section, we compute the corank of our Poisson bracket on the
space of twisted polygons.
The corank of a Poisson bracket on a smooth manifold is the codimension of the generic
symplectic leaves.
These symplectic leaves can be locally described as levels $F_i=\mathrm{const}$ of
the Casimir functions. See \cite{Wei} for the details.

For us, the only genericity condition we need is
\begin{equation}
\label{GenEq}
x_i\not=0,
\qquad
y_j\not=0; \hskip 40 pt \forall\ i,j.
\end{equation}

Our next result refers to Equation \ref{casimir}.

\begin{lemma}
\label{CasPro}
The Poisson bracket has corank 2 if $n$ is odd and corank 4 if $n$ is
even.
\end{lemma}

\proof
The Poisson bracket is quadratic in coordinates $(x,y)$.
It is very easy to see that in so-called logarithmic coordinates 
$$
p_i=\log{}x_i,
\qquad
q_i=\log{}y_i
$$
the bracket is given by a constant skew-symmetric matrix.
More precisely, the bracket between the $p$-coordinates is given by
the marix
$$
\left(
\begin{matrix}
0&-1&0&\ldots&1\\
1&0&-1&\ldots&0\\
0&1&0&\ldots&0\\
\ldots&&&\ldots&\\
-1&0&0&\ldots&0\\
\end{matrix}
\right)
$$
whose rank is $n-1$, if $n$ is odd and $n-2$, if $n$ is even.
The bracket between the $q$-coordinates is given by the opposite matrix.
\proofend

The following corollary is immediate from the preceding result
and Lemma \ref{cas}.

\begin{corollary}
\label{symL}
If $n$ is odd, then the Casimir functions are of the form $F \left(O_n,E_n \right)$.
If $n$ is even, then the Casimir functions are of the form
$F(O_{n/2},E_{n/2},O_n,E_n)$.  In both cases
the generic symplectic leaves of the Poisson structure
have dimension $4[(n-1)/2]$.
\end{corollary}

\begin{remark}
{\rm Computing the gradients, we see that a level set of the
Casimir functions is smooth as long as all the
functions are nonzero.  Thus, the generic level sets are
smooth in quite a strong sense.}
\end{remark}

\subsection{The end of the proof}

In this section, we finish the proof of Theorem 1.
Let us summarize the situation.  First we consider the
case when $n$ is odd.   On the space $\cP_n$ we have a 
generically defined and $T$-invariant
Poisson bracked that is invariant under
the pentagram map.  This bracket has co-rank $2$, and the
generic level set of the Casimir functions has dimension
$4[n/2]=2n-2$.  On the other hand, after we
exclude the two Casimirs, we have
$2[n/2]=n-1$ algebraically independent invariants that
Poisson commute with each other.  This gives us
the classical Arnold-Liouville complete integrability.

In the even case, our symplectic leaves have dimension
$4[(n-1)/2]=2n-4$.   The invariants $E_{n/2}$ and $O_{n/2}$
are also Casimirs in this case.  Once we exclude these,
we have $2[(n-1)/2]=n-2$ algebraically independent invariants.
Thus, we get the same complete integrability as in the odd
case.

This completes the proof of our Main Theorem.  In the
next chapter, we consider geometric situations where
the Main Theorem leads to quasi-periodic dynamics of
the pentagram map.

\section{Quasi-periodic motion}

In this chapter, we explain some geometric situations where our Main Theorem,
an essentially algebraic result, translates into quasi-periodic motion for the
dynamics.  The universally convex polygons furnish our main example.

\subsection{Universally convex polygons}

In this section, we define universally convex polygons and prove
some basic results about them.

We say that a matrix $M \in SL_3(\R)$ is {\it strongly diagonalizable\/} if
it has $3$ distinct positive real eigenvalues.  Such a matrix represents
a projective transformation of $\RP^2$.  We also let $M$ denote the
action on $\RP^2$.   Acting on $\RP^2$, the map $M$ fixes
$3$ distinct points.  These points, corresponding to the eigenvectors,
are in general position.  $M$ stabilizes the $3$ lines determined
by these points, taken in pairs. 
The complement of the $3$ lines is a union of $4$
open triangles.  Each open triangle is preserved by the projective
action.  We call these triangles the $M$-{\it triangles\/}.

Let $\phi \in \cP_n$ be a twisted $n$-gon, with monodromy $M$.  We
call $\phi$ {\it universally convex\/} if
\begin{itemize}
\item $M$ is a strongly diagonalizable matrix.
\item $\phi(\Z)$ is contained in one of the $M$-triangles.
\item The polygonal arc obtained by connecting consecutive vertices of
$\phi(\Z)$ is convex.
\end{itemize}
The third condition requires more explanation.  In $\RP^2$ there are
two ways to connect points by line segments.  We require the connection
to take place entirely inside the $M$-triangle that contains $\phi(\Z)$.
This determines the method of connection uniquely.

We normalize so that $M$ preserves the line at infinity
and fixes the origin in $\R^2$.   We further
normalize so that the action on $\R^2$ is given by a diagonal
matrix with eigenvalues $0<a<1<b$.  This $2 \times 2$ diagonal
matrix determines $M$.  For convenience, we will usually
work with this auxilliary $2 \times 2$ matrix.  We slightly
abuse our notation, and also refer to this $2 \times 2$ matrix
as $M$. With our normalization, the
$M$-triangles are the open quadrants in 
$\R^2$.  Finally, we normalize so that $\phi(\Z)$ is contained in
the positive open quadrant.

\begin{lemma}
\label{universal}
${\cal U\/}_n$ is open in
$\cP_n$.
\end{lemma}

\proof  
Let $\phi$ be a universally convex $n$-gon and let
$\phi'$ be a small perturbation.  Let $M'$ be
the monodromy of $\phi'$.  If the perturbation
is small, then $M'$ remains strongly diagonalizable.
We can conjugate so that $M'$ is normalized exactly
as we have normalized $M$.   

If the perturbation is small, the first $n$ points
of $\phi'(\Z)$ remain in the open positive
quadrant, by continuity.  But then all points of
$\phi'(\Z)$ remain in the open positive quadrant,
by symmetry.  This is to say that $\phi'(\Z)$ is
contained in an $M'$-triangle.

If the perturbation is small, then $\phi'(\Z)$ is
locally convex at some collection of $n$ consecutive
vertices. But then $\phi'(\Z)$ is a locally convex
polygon, by symmetry. 
The only way that $\phi(\Z)$ could fail to be convex
is that it wraps around on itself. But, the
invariance under the $2\times 2$ hyperbolic
matrix precludes this possibility.
Hence $\phi(\Z)$ is convex.
\proofend

\begin{lemma}
${\cal U\/}_n$ is invariant under the pentagram map.
\end{lemma}

\proof
Applying the pentagram map to $\phi(\Z)$ all at once, we see
that the image is again strictly convex and has the same
monodromy.
\proofend

\subsection{The Hilbert perimeter}

In this section, we introduce an invariant we call
{\it the Hilbert perimeter\/}.  This invariant plays
a useful role in our proof, given in the next
section, that the
level sets of the monodromy functions in
${\cal U\/}_n$ are compact. 

 As a prelude to
our proof, we introduce another projective
invariant -- a function of the Casimirs -- which we call the
{\it Hilbert Perimeter\/}.  This invariant is also considered
in \cite{Sch1}, and for similar purposes.

\begin{figure}[hbtp]
\centering
\includegraphics[height=1.7in]{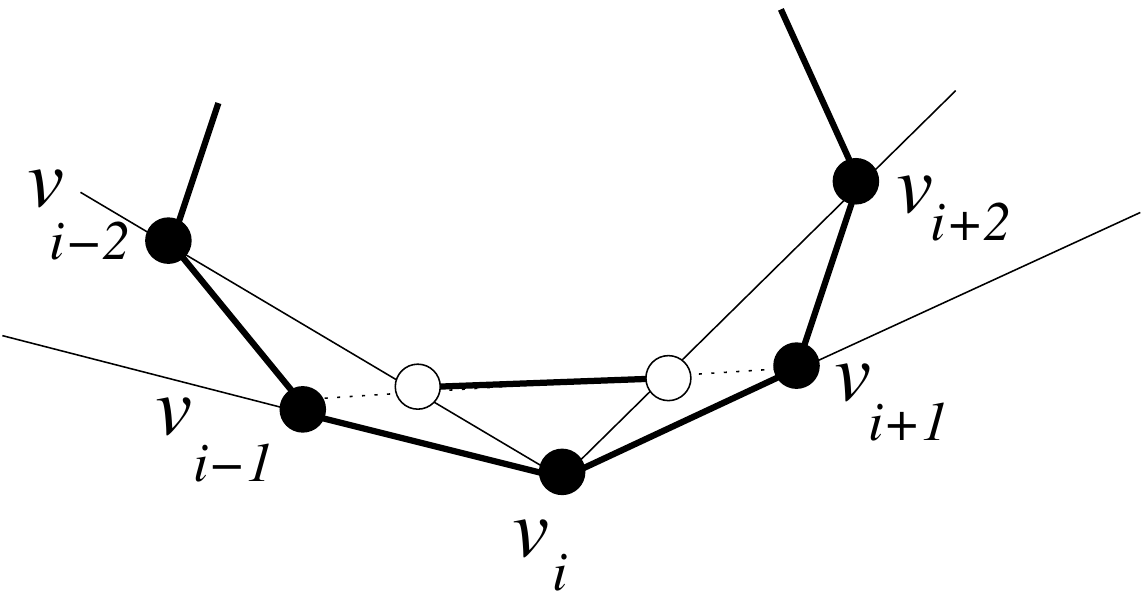}
\newline
\caption{Points involved in the definition of the invariants}
\label{ptinv}
\end{figure}

Referring to Figure \ref{corners}, we define
\begin{equation}
z_k=[(v_i,v_{i-2}),(v_i,v_{i-1}),(v_i,v_{i+1}),(v_i,v_{i+2})].
\end{equation}
We are taking the cross ratio of the slopes of the
$4$ lines in Figure \ref{ptinv}.

We now define a ``new'' invariant
\begin{equation}
H=\frac{1}{\prod_{i=1}^n z_i}
\end{equation}

\begin{remark}
{\rm Some readers will know that one can put a canonical metric
inside any convex shape, called the {\it Hilbert metric\/}.
In case $\phi$ is a genuine convex polygon, the quantity
$-\log(z_k)$ measures the Hilbert length of the thick line
segment in Figure \ref{ptinv}.  (The reader who does not know what the
Hilbert metric is can take this as a definition.)  Then
$\log(H)$ is the Hilbert perimeter of $T(P)$ with respect to
the Hilbert metric on $P$.  Hence the name.
}
\end{remark}

\begin{lemma}
$H=1/(O_nE_n)$.  
\end{lemma}

\proof
This is a local calculation, which amounts to showing that
$z_k=x_ky_k$. The best way to do the calculation is to normalize
so that $4$ of the points are the vertices of a square.  We
omit the details.
\proofend

\subsection{Compactness of the level sets}

In this section, we prove that the level sets of the
monodromy functions in ${\cal U\/}$ are compact.

Let ${\cal U\/}_n(M,H)$ denote the subset of ${\cal U\/}_n$
consisting of elements whose monodromy is $M$ and whose
Hilbert Perimeter is $H$.  
In this section we will prove
that ${\cal U\/}_n(M,H)$ is compact.  
For ease of notation, we abbreviate this space by $X$. 

Let $\phi \in X$.
We normalize so that 
$M$ is as in Lemma \ref{universal}.   We also normalize so that
$\phi(0)=(1,1)$.  Then there are numbers $(x,y)$ such that
$\phi(n)=(x,y)$, independent of the choice of $\phi$.
We can assume that $x>1$ and $y<1$.   The portion of
$\phi$ of interest to us, namely
$\phi(\{0,...,n-1\})$, lies entirely in the rectangle $R$
whose two opposite corners are $(1,1)$ and $(x,y)$.
Let $(v_i,v_j)$ denote the line determined by
$v_i$ and $v_j$.  Here $v_k=\phi(k)$.
In particular, let $L_i=(v_i,v_{i+1})$.

\begin{lemma}
Suppose that $\{\phi_k\} \in X$  is a sequence that
does not converge on a subsequence to another element of $X$.
Then,  passing to a subsequence we can arrange that
at least one of the two situations holds: there exists some $i$ such that
\begin{itemize}
\item The angle between $L_i$ and $L_{i+1}$ tends to $0$ as
$k \to\infty$ whereas
the angle between $L_{i+1}$ and $L_{i+2}$ does not;
\item The points $v_i$ and $v_{i+1}$ converge to a common point
as $k \to \infty$ whereas
$v_{i+2}$ converges to a distinct point.
\end{itemize}
\end{lemma}

\proof
Suppose that there is some minimum distance $\epsilon$ between
all points of $\phi_k$ in the rectangle $R$.
In this case, the angle between two consecutive segments
must tend to $0$ as $k \to \infty$.  However, not all angles between
consecutive segments can converge to $0$ because of the
fixed monodromy.  The first case is now easy to arrange.
If there is no such minimum $\epsilon$, then two points
coalesce, on a subsequence.  For the same reason as above,
not all points can coalesce to the same point.  The second
case is now easy to arrange.
\proofend

\begin{lemma} $X$ is compact.
\end{lemma}

We will suppose we have the kind of sequence we had in the previous
lemma and then  derive a contradiction.  In the first case above,
the slopes of the lines $(v_{i+2},v_{i})$ and
$(v_{i+2},v_{i+1})$ converge to each other as $k \to \infty$,
but the common limit remains uniformly bounded away from
the slopes of $(v_{i+2},v_{i+3})$ and $(v_{i+2},v_{i+4})$. Hence
$z_{i+2} \to 0$.  Since $z_j \in (0,1)$ for all $j$, we have
$H \to \infty$ in this case.  This is a contradiction.

To deal with the second case, we can assume that the first
case cannot be arranged.  That is, we can assume that there
is a uniform lower bound to the angles between two consecutive
lines $L_i$ and $L_{i+1}$ for all indices and all $k$.
But then the same situation as in Case 1 holds, and we get
the same contradiction.
\proofend

\subsection{Proof of Theorem \ref{universally convex}}

In this section, we finish the proof of Theorem \ref{universally convex}.

Recall that the level sets of our Casimir functions give
a (singular) foliation by symplectic leaves.
Note that all corner invariants are nonzero for points in
${\cal U\/}_n$.  Hence, our singular symplectic foliation
intersects ${\cal U\/}_n$
 in leaves that are all smooth symplectic manifolds.
Let $k=[(n-1)/2]$.  

Let $\cal M$ be a symplectic leaf.  Note
that $\cal M$ has dimension $4k$. 
Consider the map
\begin{equation}
F=(O_1,E_1,...,O_k,E_k),
\end{equation}
made from our algebraically independent monodromy
invariants.  Here we are excluding all the Casimirs
from the definition of $F$.

Say that a point $p \in \cal M$ is {\it regular\/} if
$dF_p$ is surjective.  Call $\cal M$ {\it typical\/}
if some point of $\cal M$ is regular.
Given our algebraic independence result, and the
fact that the coordinates of $F$ are polynomials, we see that
almost every  symplectic leaf is typical.

\begin{lemma}
If $\cal M$ is typical then almost every $F$-fiber of $\cal M$
is a smooth submanifold of $\cal M$.
\end{lemma}

\proof
Let $S=F({\cal M\/}) \subset \R^{2k}$.  Note that
$S$ has positive measure since $dF_p$ is nonsingular
for some $p \in {\cal M\/}$.  Let $\Sigma \subset \cal M$
denote the set of points $p$ such that $dF_p$ is
not surjective. Sard's theorem says that
$F(\Sigma)$ has measure $0$.   Hence, almost every
fiber of $\cal M$ is disjoint from $\Sigma$.
\proofend

Let $\cal M$ be a typical symplectic leaf, and let
$\cal F$ be a smooth fiber of $F$.  Then $\cal F$
has dimension $2k$. Combining our Main Theorem with
the standard facts about Arnold-Liouville complete
integrability (e.g., \cite{Arn}), we see that the
monodromy invariants give a canonical affine
structure to $\cal F$.   The pentagram map $T$
preserves both $\cal F$, and is a translation
relative to this affine structure.  Any pre-compact
orbit in $\cal F$ exhibits quasi-periodic motion.

Now,
$T$ also preserves the monodromy.  But then each $T$-orbit
in $\cal F$ is contained in one of our spaces
${\cal U\/}_n(H,M)$.  Hence, the orbit is precompact.
Hence, the orbit undergoes quasi-periodic motion.
Since this argument works for almost every $F$-fiber
of almost every symplectic leaf in ${\cal U\/}_n$,
we see that almost every orbit in ${\cal U\/}_n$
undergoes quasi-periodic motion under the pentagram map.

This completes the proof of Theorem \ref{universallyconvex}.

\begin{remark}
{\rm We can say a bit more.
For almost every choice of monodromy $M$,
the intersection
\begin{equation}
{\cal F\/}(M)={\cal F\/} \cap {\cal U\/}_n(H,M)
\end{equation}
is a smooth compact submanifold and inherits an
invariant affine structure from $\cal F$.  In
this situation, the restriction of $T$ to 
${\cal F\/}(M)$ is a translation in the
affine structure.}
\end{remark}

\subsection{Hyperbolic cylinders and tight polygons}

In this section, we put Theorem \ref{universallyconvex} in a
somewhat broader context.  The material in this section
is a prelude to our proof, given in the next section,
of a variant of Theorem \ref{universallyconvex}.

Before we sketch variants of Theorem \ref{universally convex}, we
think about these polygons in a different way.
A {\it projective cylinder\/} is a topological cylinder that has
coordinate charts into $\RP^2$ such that the transition functions
are restrictions of projective transformations.  This is a classical
example of a {\it geometric structure\/}. See \cite{thu} or \cite{OT} for details.

\begin{example}
{\rm Suppose that $M$ acts on $\R^2$ as a nontrivial
diagonal matrix having eigenvalues $0<a<1<b$.  Let $Q$ denote the
open positive quadrant.  Then $Q/M$ is a projective cylinder.
We call $Q/M$ a {\it hyperbolic cylinder\/}.}
\end{example}

Let $Q/M$ be a hyperbolic cylinder.  Call a polygon on
$Q/M$ {\it tight\/} if it has the following $3$ properties.
\begin{itemize}
\item It is embedded;
\item It is locally convex;
\item It is homologically nontrivial.
\end{itemize}
Any universally convex polygon gives rise to a tight polygon on
$Q/M$, where $M$ is the monodromy normalized in the standard way.  The
converse is also true.   Moreover, two tight polygons on $Q/M$ give
rise to equivalent universally convex polygons iff some
locally projective diffeomorphism of $Q/M$ carries one
polygon to the other.  We call such maps {\it automorphisms\/}
of the cylinder, for short.

 Thus, we can think of the pentagram map
as giving an iteration on the space of tight polygons on a
hyperbolic cylinder.  There are $3$ properties that give rise
to our result about periodic motion.
\begin{enumerate}
\item The image of a tight polygon under the pentagram map is another
well-defined tight polygon.
\item The space of tight polygons on a hyperbolic cylinder, modulo the
projective automorphism group, is compact.
\item The strongly diagonalizable elements are open in $SL_3(\R)$.
\end{enumerate}
The third condition guarantees that the set of all
tight polygons on all hyperbolic cylinders is an open
subset of the set of all twisted polygons.

\subsection{A related theorem}

In this section we prove a variant of
Theorem \ref{universally convex} for a different
family of twisted polygons.

\begin{figure}[hbtp]
\centering
\includegraphics[height=1.3in]{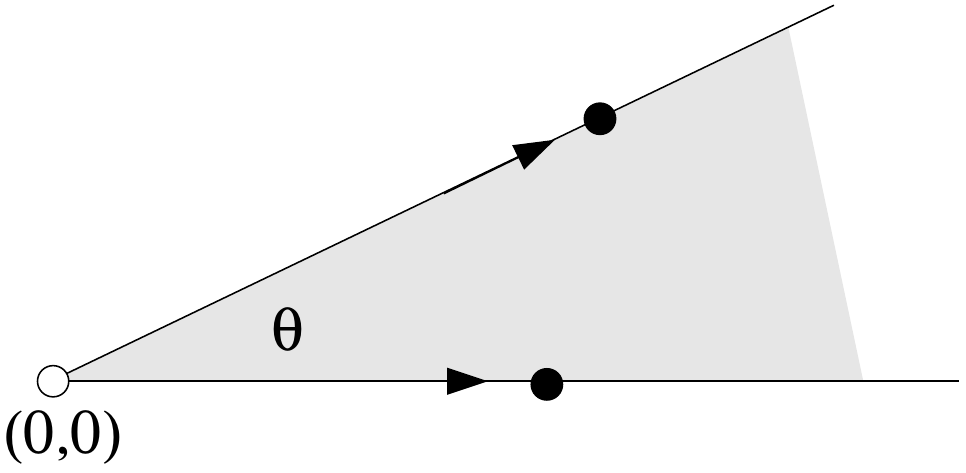}
\newline
\caption{The cylinder $\Sigma(\theta,d)$}
\label{cylinder}
\end{figure}

We start with a sector of angle $\theta$ in the plane, as shown in Figure
\ref{cylinder}, and glue the top edge to the bottom edge by a similarity
$S$ that has dilation factor $d$.   We omit the origin
from the sector.  The quotient is the projective cylinder
we call $\Sigma(\theta,d)$.  When $d=1$ we have a Euclidean
cone surface.  When $\theta=2 \pi$ we have the punctured plane.

We consider the case when $\theta$ is small and $d$ is close to $1$.
In this case, $\Sigma(\theta,d)$ admits tight polygons for any
$n$.  (It is easiest to think about the case when $n$ is large.)
When developed out in the plane, these tight polygons follow along
logarithmic spirals.

Let $S(\theta,d)$ denote the subset of $\R^2$ consisting of pairs
$(\theta',d')$ where
\begin{equation}
0<\theta'<\theta; \hskip 30 pt
1<d'<d
\end{equation}

Define
\begin{equation}
\widehat \Sigma(d,1)=\bigcup_{(\theta',d') \in S(\theta,d)}\Sigma(\theta',d')
\end{equation}
One might say that $\widehat \Sigma(\theta,d)$ is the space of polygons
that are more tightly coiled than those on $\Sigma(\theta,d)$.

\begin{theorem}
\label{variant}
Suppose that $\theta>0$ is sufficiently close to $0$ and $d>1$ is
sufficiently close to $1$.  Then almost every point of
$\widehat \Sigma(\theta,d)$ lies on a smooth torus that has a
$T$-invariant affine structure.  Hence, the orbit of almost
every point of $\widehat \Sigma(\theta,d)$ undergoes
quasi-periodic motion.
\end{theorem}

\proof
Our proof amounts to verifying the three properties above for the
points in our space.  We fix $(\theta,d)$ and let
$\widehat \Sigma=\widehat \Sigma(\theta,d)$. 
\begin{enumerate}
\item
Let $P$ be a tight polygon on $\Sigma(\theta',d')$.
If $\theta'$ is sufficiently small and $d'$ is sufficiently close to $1$, then
each vertex $v$ of $P$ is much closer to its
neighbors than it is to the origin.  For this reason, the pentagram
map acts on, and preserves, the set of tight polygons on $\Sigma(\theta',d')$.
The same goes for the inverse of the pentagram map.  Hence
$\widehat \Sigma$ is a $T$-invariant subset of
$\cP_n$.
\item 
Let $Z(\theta',d',\alpha)$ denote the space of
tight polygons on $\Sigma(\theta',d')$ having
Hilbert perimeter $\alpha$.  We consider these
tight polygons equivalent if there is a similarity of
$\Sigma(\theta',d')$ that carries one to the other.
A proof very much like the compactness argument given in
\cite{Sch1}, for closed polygons, shows that
$Z(\theta',d',\alpha)$ is compact for
$\theta'$ near $0$ and $d'$ near $1$ and $\alpha$ arbitrary.
Hence, the level sets 
of the Casimir functions intersect $\widehat \Sigma$
in compact sets.
\item 
The similarity $S$ is the monodromy for our tight
polygons. $S$ lifts to an element of $SL_3(\R)$
that has one real eigenvalue and two complex conjugate eigenvalues.
Small perturbations of $S$ have the same property.  Hence,
$\widehat \Sigma$ is open in
$\cP_n$.
\end{enumerate}
We have assembled all the ingredients necessary for the
proof of Theorem \ref{universally convex}.  The same argument
as above now establishes the result.
\proofend
\begin{remark}
{\rm The first property crucially uses the fact that
$\theta$ is small.  Consider the case $\theta=2 \pi$. It
can certainly happen that $P$ contains the origin in
its hull but $T(P)$ does not.
We do not know the exact bounds on $\theta$ and $d$ necessary
for this construction.}
\end{remark}

\section{Another coordinate system in space $\cP_n$} \label{another}

\subsection{Polygons and difference equations}\label{DiffEq}

Consider two arbitrary $n$-periodic sequences $(a_i),\,(b_i)$ with  $a_i,b_i\in\R$
and $i\in\Z$, such that $a_{i+n}=a_i,\,b_{i+n}=b_i$.
Assume that $n\not=3\,m$. This will be our standing assumption whenever we work with the $(a,b)$-coordinates; its meaning will become clear shortly.
We shall associate to these sequences a difference
equation of the form
\begin{equation} 
\label{recur}
V_{i+3}=a_i \,V_{i+2}+b_i\,V_{i+1} + V_i,
\end{equation}
for all $i$.

A solution $V=(V_i)$ is a sequence of numbers $V_i\in\R$ satisfying (\ref{recur}).
Recall a well-known fact that
the space of solutions of (\ref{recur}) is 3-dimensional (any solution is determined by the
initial conditions $(V_0,V_1,V_2)$).
We will often understand $V_i$ as vectors in $\R^3$.
The $n$-periodicity then implies that there exists a matrix $M\in\SL(3,\R)$
called the monodromy matrix, such that $V_{i+n}=M\,V_i$.

\begin{proposition} 
\label{basic}
If $n$ is not divisible by 3 then the space $\cP_n$ is isomorphic to
the space of the equations (\ref{recur}).
\end{proposition}

\proof
First note that since $\PGL(3,\R)\cong\SL(3,\R)$,
every $M\in\PGL(3,\R)$ corresponds to a unique  element of $\SL(3,\R)$
that (abusing the notations) we also denote by $M$.

{\bf A}.
Let $(v_i),i\in\Z$ be a sequence of points $v_i\in\RP^2$ in general position
with monodromy~$M$.
Consider first an arbitrary lift of the points $v_i$ to vectors $\tilde V_i\in\R^3$
with the condition $\tilde V_{i+n}=M(\tilde V_i)$.
The general position property implies that 
$\det(\tilde V_i,\tilde V_{i+1},\tilde V_{i+2})\not=0$ for all $i$.
The vector $\tilde V_{i+3}$ is then a linear combination of the linearly independent vectors 
$\tilde V_{i+2}, \tilde V_{i+1}, \tilde V_i$, that is,
$$
\tilde V_{i+3}=a_i \,\tilde V_{i+2}+b_i \,\tilde V_{i+1} + c_i\,\tilde V_i,
$$
for some $n$-periodic sequences $(a_i),\,(b_i),\,(c_i)$.
We wish to rescale: $V_i=t_i \tilde V_i$, so that
\begin{equation} 
\label{Wron}
\det(V_i,V_{i+1},V_{i+2})=1
\end{equation}
for all $i$.
Condition (\ref{Wron}) is equivalent to $c_i\equiv1$.
One obtains the following system of equations in $(t_1,\ldots,t_n)$:
$$
\begin{array}{rcll}
t_it_{i+1}t_{i+2}&=&1/\det(\tilde V_i,\tilde V_{i+1},\tilde V_{i+2}),
&
i=1,\dots,n-2\\[5pt]
t_{n-1}t_{n}t_{1}&=&
1/\det(\tilde V_{n-1},\tilde V_{n},\tilde V_{1}),&\\[5pt]
t_{n}t_{1}t_{2}&=&
1/\det(\tilde V_{n},\tilde V_{1},\tilde V_{2})&
\end{array}
$$
This system has a unique solution if $n$ is not divisible by 3.
This means that any generic twisted $n$-gon in $\RP^2$ has a unique lift to
$\R^3$ satisfying (\ref{Wron}).
We proved that a twisted $n$-gon  defines an equation (\ref{recur})
with $n$-periodic $a_i,b_i$.

Furthermore, if $(v_i)$ and $(v'_i),i\in\Z$ are two 
projectively equivalent twisted  $n$-gons, then
they correspond to the same equation (\ref{recur}).
Indeed, there exists $A\in\SL(3,\R)$ such that
$A(v_i)=v'_i$ for all $i$.
One has, for the (unique) lift: $V'_i=A(V_i)$.
The sequence $(V'_i)$ then obviously satisfies the same equation
(\ref{recur}) as $(V_i)$.

{\bf B}.
Conversely, let $(V_i)$ be a sequence of vectors $V_i\in\R^3$ satisfying (\ref{recur}).
Then every three consecutive points satisfy (\ref{Wron}) and, in particular,
are linearly independent.
Therefore, the projection $(v_i)$ to $\RP^2$ satisfies the general position condition.
Moreover, since the sequences $(a_i),\,(b_i)$ are $n$-periodic,
$(v_i)$ satisfies $v_{i+n}=M(v_i)$.
It follows that every equation (\ref{recur}) defines a generic twisted $n$-gon.
A choice of initial conditions $(V_0,V_1,V_2)$ fixes a twisted polygon, a different choice yields a projectively equivalent one.
\proofend

Proposition \ref{basic} readily implies the next result.

\begin{corollary}
\label{Space}
If $n$ is not divisible by 3 then $\cP_n =\R^{2n}$.
\end{corollary}

We call the lift $(V_i)$ of the sequence $(v_i)$ satisfying equation (\ref{recur}) with $n$-periodic $(a_i,b_i)$ 
{\it canonical}.

\begin{remark}
{\rm
The isomorphism between the space $\cP_n$ and the space of difference equations (\ref{recur}) (for $n\neq 3m$) goes back to the classical ideas of projective differential geometry.
This is a discrete version of the well-known isomorphism 
between the space of smooth non-degenerate curves in $\RP^2$
and the space of linear differential equations, see \cite{OT} and references therein
and Section \ref{Curva}. The ``arithmetic restriction'' $n\not=3m$ is quite remarkable.

Equations (\ref{recur}) and their analogs were already used in \cite{FRS}
in the context of integrable systems;
in the $\RP^1$-case these equations were recently considered in \cite{Mar}
to study the discrete versions of the Korteweg - de Vries equation. It is notable that an analogous
arithmetic assumption $n\neq 2m$ is made in this paper as well.
}
\end{remark}

\begin{remark}
{\rm Let us now comment on what happens if $n$ is divisible by 3. A certain modification of Proposition \ref{basic} holds in this case as well. Given a twisted $n$-gon $(v_i)$ with monodromy $M$, lift points $v_0$ and $v_1$ arbitrarily as vectors $V_0,V_1 \in \R^3$, and then continue lifting consecutive points so that the determinant condition (\ref{Wron}) holds. This implies that equation (\ref{recur}) holds as well. 

One has: 
\begin{equation} \label{MV}
M(V_i)=t_i V_{i+n}
\end{equation}
 for non-zero reals $t_i$, and (\ref{Wron}) implies that $t_it_{i+1}t_{i+2}=1$ for all $i\in \Z$. It follows  that the sequence $t_i$ is 3-periodic; let us write $t_{1+3j}=\alpha, t_{2+3j}=\beta, t_{3j}=1/(\alpha\beta)$. Applying the monodromy linear map $M$ to (\ref{recur}) and using (\ref{MV}), we conclude that 
$$
a_{n+i}=\frac{t_{i+2}}{t_i} a_i, \ \ b_{n+i}=\frac{t_{i+1}}{t_i} b_i,
$$
that is, 
\begin{equation} \label{scalea}
\begin{split}
a_{n+3j}= \alpha\beta^2\ a_{3j}, \ a_{n+3j+1}=\frac{1}{\alpha^2\beta}\ a_{3j+1}, \ a_{n+3j+2}=\frac{\alpha}{\beta}\ a_{3j+2}, \\
b_{n+3j}= \alpha^2\beta \ b_{3j}, \ b_{n+3j+1}=\frac{\beta}{\alpha} \ b_{3j+1}, \ b_{n+3j+2}=\frac{1}{\alpha\beta^2}\ b_{3j+2}.
\end{split}
\end{equation}
We are still free to rescale  $V_0$ and $V_1$. This defines an action of the group $\R^* \times \R^*$:
$$
V_0 \mapsto u V_0,\ V_1 \mapsto v V_1,\ \ u\neq 0, v\neq 0. 
$$
The action of the group $\R^* \times \R^*$ on the coefficients $(a_i,b_i)$ is as follows:
\begin{equation} \label{action}
\begin{split}
a_{3j} \mapsto u^2v\ a_{3j}, \ a_{3j+1} \mapsto \frac{v}{u}\ a_{3j+1},\ a_{3j+2} \mapsto \frac{1}{uv^2}\ a_{3j+2},\\
b_{3j} \mapsto \frac{u}{v}\ b_{3j}, \ b_{3j+1} \mapsto uv^2\ b_{3j+1},\ b_{3j+2} \mapsto \frac{1}{u^2v}\ b_{3j+2}.
\end{split}
\end{equation}
When  $n\neq 3m$, according to (\ref{scalea}), this action makes it possible to normalize all $t_i$ to $1$ which makes the lift canonical. However, if $n=3m$ then the $\R^* \times \R^*$-action on $t_i$ is trivial, and the pair $(\alpha,\beta)\in \R^* \times \R^*$  is a projective
 invariant of the twisted polygon. One concludes that $\cP_n$ is the orbit space
$$
[\{(a_0,\dots,a_{n-1},b_0,\dots,b_{n-1})\}/(\R^* \times \R^*)] \times (\R^* \times \R^*)
$$
with respect to  $\R^* \times \R^*$-action (\ref{action}). This statement replaces Proposition \ref{basic} in the case of $n=3m$.
 
It would be interesting to understand the geometric meaning of the ``obstruction" $(\alpha,\beta)$. If the obstruction is trivial, that is, if $\alpha=\beta=1$, then there exists a 2-parameter family of canonical lifts, but if the obstruction is non-trivial then no canonical lift exists.
}
\end{remark}

\subsection{Relation between the two coordinate systems}

We now have two coordinate systems, $(x_i,y_i)$ and $(a_i,b_i)$. Assuming that $n$ is not divisible by 3, let us calculate the relations between the two systems.

\begin{lemma}
\label{atox}
One has:
\begin{equation}
\label{RelEq}
x_{i}=\frac{a_{i-2}}{b_{i-2}\,b_{i-1}}, 
\qquad
y_{i}=-\frac{b_{i-1}}{a_{i-2}\,a_{i-1}}.
\end{equation}
\end{lemma}

\proof
Given four vectors $a,b,c,d$ in $\R^3$, the intersection line of the planes 
$\mathrm{Span}(a,b)$ and $\mathrm{Span}(c,d)$ is spanned 
by the vector $(a\times b)\times (c\times d)$.
Note that the volume element equipes $\R^3$ with the bilinear vector product:
$$
\R^3\times{}\R^3\to\left(\R^3\right)^\star.
$$
Using the identity 
\begin{equation} 
\label{ident}
(a\times b)\times (b\times c)=\det(a,b,c)\,b,
\end{equation}
and the recurrence (\ref{recur}), 
let us compute lifts of the quadruple of points 
$$
\left(
v_{i-1},\,v_i,\,(v_{i-1},v_i)\cap(v_{i+1},v_{i+2}),\,(v_{i-1},v_i)\cap(v_{i+2},v_{i+3})
\right)
$$
involved in the left corner cross-ratio.
One has
$$
V_{i-1}=V_{i+2}-a_{i-1}\,V_{i+1}-b_{i-1}\,V_i.
$$
Furthermore, it is easy to obtain the lift of the intersection points 
involved in the left corner cross-ratio.
For instance, $(v_{i-1},v_i)\cap(v_{i+1},v_{i+2})$ is
$$
\begin{array}{rcl}
\left(
V_{i-1}\times{}V_i
\right)
\times
\left(
V_{i+1}\times{}V_{i+2}
\right)
&=&
\left(
(V_{i+2}-a_{i-1}\,V_{i+1}-b_{i-1}\,V_i)\times{}V_i
\right)
\times
\left(
V_{i+1}\times{}V_{i+2}
\right)\\[6pt]
&=&
V_{i+2}-a_{i-1}\,V_{i+1}.
\end{array}
$$
One finally obtains the following four vectors in $\R^3$:
$$
\left(
V_{i+2}-a_{i-1}\,V_{i+1}-b_{i-1},V_i, 
\quad
V_i, 
\quad
V_{i+2}-a_{i-1}\,V_{i+1}, 
\quad
b_i\,V_{i+2}-a_{i-1}\,V_i-a_{i-1}\,b_i\,V_{i+1}
\right).
$$
Similarly, for the points involved in the right corner cross-ratio
$$
\left(
 a_i\,V_{i+2}+b_i,V_{i+1}+V_i,
 \quad
 V_{i+2},
 \quad
 b_i\,V_{i+1}+V_i,
 \quad
 b_i\,V_{i+2}-a_{i-1}\,V_i-a_{i-1}\,b_i\,V_{i+1}
 \right).
$$

Next, given four coplanar vectors $a,b,c,d$ in $\R^3$ such that
$$
c=\lambda_1\,a+\lambda_2\,b,
\qquad
d=\mu_1\,a+\mu_2\,b,
$$
where $\lambda_1,\lambda_2,\mu_1,\mu_2$ are arbitrary constants,
the cross-ratio of the lines spanned by these vectors is given by
$$
[a,b,c,d]=\frac{\lambda_2\mu_1-\lambda_1\mu_2}{\lambda_2 \mu_1}.
$$
Applying this formula to the two corner cross-ratios yields the result.
\proofend

Formula (\ref{RelEq}) implies the following relations:
\begin{equation}
\label{xtoa}
x_{i}\,y_{i}=-\frac{1}{a_{i-1}\,b_{i-2}}, 
\quad
x_{i+1}\,y_{i}=-\frac{1}{a_{i-2}\,b_{i}},
\quad
\frac{a_{i}}{a_{i-3}}=\frac{x_{i}\,y_{i-1}}{x_{i+1}\,y_{i+1}},
\quad
\frac{b_{i}}{b_{i-3}}=\frac{x_{i-1}\,y_{i-1}}{x_{i+1}\,y_{i}}
\end{equation}
that will be of use later.

\begin{remark} \label{sense}
{\rm If $n$ is a multiple of 3 then the coefficients $a_i$ and $b_i$ are not well defined and they are not $n$-periodic anymore; however, according to formulas (\ref{scalea}) and (\ref{action}), the right hand sides of formulas (\ref{RelEq}) are still well defined and are $n$-periodic.
}
\end{remark}

\subsection{Two versions of the projective duality}

We now wish to express the pentagram map $T$ in the $(a,b)$-coordinates. 
We shall see that  $T$ is the composition of two involutions each of which is a kind of projective duality.

The notion of projective duality in $\RP^2$ is based on the fact that the dual projective plane
$\left(\RP^2\right)^\star$ is the space of one-dimensional subspaces of $\RP^2$
which is again equivalent to $\RP^2$.
Projective duality applies to smooth curves: it associates to a curve
$\g(t)\subset\RP^2$ the 1-parameter family of its tangent lines.
In the discrete case, there are different ways to define projectively dual polygons.
We choose  two simple versions.

\begin{figure}[hbtp]
\centering
\includegraphics[width=5in]{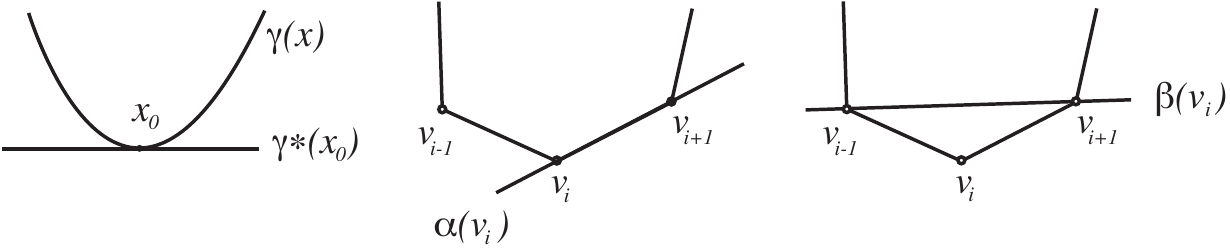}
\newline
\caption{Projective dual for smooth curves and polygons}
\label{DualFig}
\end{figure}

\begin{definition}
{\rm
Given a sequence of points $v_i \in \RP^2$, we define two sequences
$\a(v_i)\in(\RP^2)^\star$ and $\b(v_i)\in(\RP^2)^\star$ as follows:

\begin{enumerate}
\item
$\a(v_i)$ is the line $(v_i,v_{i+1})$,
\item
$\b(v_i)$ is the line $(v_{i-1},v_{i+1})$,
\end{enumerate}

\noindent
see Figure \ref{DualFig}.
}
\end{definition}

Clearly, $\a$ and $\b$ commute with the natural $\PGL(3,\R)$-action and therefore
are well-defined on the space $\cP_n$.
The composition of $\a$ and $\b$ is precisely 
the pentagram map $T$.

\begin{lemma}
\label{TmapDual}
One has
\begin{equation}
\label{InvEq}
\a^2=\tau,
\qquad
\b^2=\mathrm{Id},
\qquad
\a\circ\b=T.
\end{equation}
where $\tau$ is the cyclic permutation:
\begin{equation}
\label{Cycle}
\tau(v_i)=v_{i+1}.
\end{equation}
\end{lemma}

\begin{figure}[hbtp]
\centering
\includegraphics[width=6in]{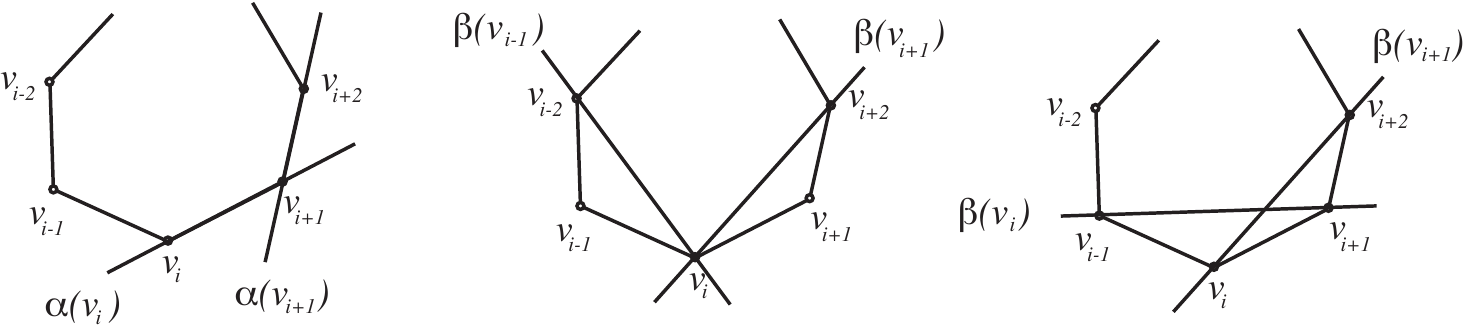}
\caption{Iteration of the duality maps:\newline
$
\a^2(v_i)=\a(v_i)\cap\a(v_{i+1}),
\quad
\b^2(v_i)=\b(v_{i-1})\cap\b(v_{i+1})
\quad\hbox{and}\quad
(\a\circ\b)(v_i)=\b(v_{i})\cap\b(v_{i+1}).
$
}
\label{IntersFig}
\end{figure}

\proof
The composition of the maps $\a$ and $\b$, with themselves and with each other, 
associates to the corresponding lines
(viewed as points of $(\RP^2)^\star$) their intersections,
see Figure \ref{IntersFig}. \proofend

The map (\ref{Cycle}) defines the natural action of the  group
$\Z$ on $\cP_n$.
All the geometric and algebraic structures we consider
 are invariant with respect to this action.

\subsection{Explicit formula for $\a$}

It is easy to calculate the explicit formula of the map $\a$ in terms of the coordinates $(a_i,b_j)$. As usual, we assume $n\not=3\,m$.

\begin{lemma}
\label{ExplAlpha}
Given a twisted  $n$-gon with monodromy $(v_i),\,i\in\Z$ represented by a difference equation (\ref{recur}),
the $n$-gon $(\a(v_i)),\,i\in\Z$ is represented by the equation (\ref{recur}) with
coefficients
\begin{equation}
 \label{recurU}
\a^*(a_i)=-b_{i+1},
\qquad
\a^*(b_i)=-a_i,
\end{equation}
where, as usual, $a^*$ stands for the pull-back of the coordinate functions.
\end{lemma}

\proof
Consider the canonical lift $(V_i)$ to $\R^3$.
Let $U_i=V_i\times V_{i+1}\in(\R^3)^\star$.
This is obviously a lift of the sequence $(\a(v_i))$ to $(\R^3)^\star$.
We claim that $(U_i)$ is, in fact, a canonical lift.

Indeed, 
$U_i$ is a lift of $u_i$ since $V_i\times V_{i+1}$ is orthogonal to $V_i$ and to $V_{i+1}$.
Next, using the identity (\ref{ident}) one has
$$
\det(U_i\times U_{i+1},U_{i+1}\times U_{i+2},U_{i+2}\times U_{i+3})=
[(U_i\times U_{i+1})\times(U_{i+1}\times U_{i+2})]\cdot (U_{i+2}\times U_{i+3})
$$
$$
=U_{i+1}\cdot (U_{i+2}\times U_{i+3})=\det(U_{i+1},U_{i+2},U_{i+3})=1.
$$
 
It follows that the sequence $U_i\in\R^3$ satisfies the equation
$$
U_{i+3}=\a^*(a_i)\,U_{i+2}+\a^*(b_i)\,U_{i+1} + U_i
$$
with some $\a^*(a_i)$ and $\a^*(b_i)$.
Let us show that these coefficients are, indeed, given by (\ref{recurU}).
For all $i$, one has 
$$
U_{i+1}\cdot V_i=1,
\qquad
U_{i}\cdot V_{i+2}=1,
\qquad
U_{i+3}\cdot V_{i+3}=0.
$$
Using  (\ref{recur}), the last identity leads to: 
$$
\a^*(b_i)\,U_{i+1}\cdot V_i +a_i\,U_{i}\cdot V_{i+2}=0.
$$
Hence $\a^*(b_i)=-a_i$. 
The first identity in (\ref{recurU}) follows from formula (\ref{InvEq}).
Indeed, one has $\a^*(\a^*(a_i))=a_{i+1}$ and $\a^*(\a^*(b_i))=b_{i+1}$, and we are done.
\proofend

\subsection{Recurrent formula for $\b$}

The explicit formula for the map $\b$ is more complicated, and
we shall give a recurrent expression.

\begin{lemma}
\label{ExplBeta}
Given an $n$-gon $(v_i),\,i\in\Z$ represented by a difference equation (\ref{recur}),
the $n$-gon $(\b(v_i)),\,i\in\Z$ is represented by the equation (\ref{recur}) with
coefficients
\begin{equation}
 \label{recurUBet}
\b^*(a_i)=-\frac{\lambda_i\,b_{i-1}}{\lambda_{i+2}},
\qquad
\b^*(b_i)=-\frac{\lambda_{i+3}\,a_{i+1}}{\lambda_{i+1}}. 
\end{equation}
where the coefficients $\lambda_i$ are uniquelly defined by
\begin{equation}
 \label{LamB}
\lambda_{i}\lambda_{i+1}\lambda_{i+2}=-\frac{1}{1+b_{i-1}a_i}
\end{equation}
for all $i$.
\end{lemma}

\proof
The lift of the map $\b$ to $\R^3$ takes $V_i$ to 
$W_i=\lambda_i V_{i-1}\times V_{i+1}$ where the coefficients 
$\lambda_i$ are chosen in such a way  that $\det(W_i,W_{i+1},W_{i+2})=1$ for all $i$.
The sequence $W_i\in\R^3$ satisfies the equation
$$
W_{i+3}=\b^*(a_i)\,W_{i+2}+\b^*(b_i)\,W_{i+1} + W_i.
$$
To find $\b^*(a_i)$ and $\b^*(b_i)$, one substitutes $W_i=\lambda_i V_{i-1}\times V_{i+1}$, and
then, using (\ref{recur}), expresses each $V$ as a linear combination of $V_i, V_{i+1}, V_{i+2}$.
The above equation is then equivalent to the following one:
$$
\begin{array}{rcl}
\left(
\b^*(a_i)\,\lambda_{i+2}+b_{i-1}\,\lambda_i
\right)
V_i\times{}V_{i+1}&&\\[6pt]
+\left(
a_{i+1}\,\lambda_{i+3}+\b^*(b_i)\,\lambda_{i+1}
\right)V_i\times{}V_{i+2}&&\\[6pt]
+\left(
(1+b_i\,a_{i+1})\,\lambda_{i+3}+\b^*(a_i)\,a_i\,\lambda_{i+2}-\lambda_i
\right)V_{i+1}\times{}V_{i+2}&=&0.
\end{array}
$$
Since the three terms are linearly independent, one obtains three relations.
The first two equations lead to (\ref{recurUBet}) while the last one
gives the recurrence 
$$
\lambda_{i+3}=
\lambda_i\,\frac{1+a_{i}\,b_{i-1}}{1+a_{i+1}\,b_i}.
$$
On the other hand, one has
$$
\lambda_{i}\,\lambda_{i+1}\,\lambda_{i+2}\,
\det\left(
V_{i-1}\times V_{i+1},V_{i}\times V_{i+2},V_{i+1}\times V_{i+3}
\right)=1.
$$
Once again, expressing each $V$ as a linear combination of $V_i, V_{i+1}, V_{i+2}$, 
yields 
$$
\lambda_{i}\lambda_{i+1}\lambda_{i+2}\,(1+a_ib_{i-1})=-1
$$
and one obtains (\ref{LamB}).
\proofend

\subsection{Formulas for the pentagram map}

We can now describe the pentagram map in terms of $(a,b)$-coordinates and to 
deduce formulas (\ref{ExpXEq}).

\begin{proposition}
\label{mapinx1}
(i)
One has:
$$
T^*(x_i)=x_i\,\frac{1-x_{i-1}\,y_{i-1}}{1-x_{i+1}\,y_{i+1}},
\qquad 
T^*(y_i)=y_{i+1}\,\frac{1-x_{i+2}\,y_{i+2}}{1-x_{i}\,y_{i}}.
$$

(ii)
Assume that $n=3m+1$ or $n=3m+2$; in both cases,
\begin{equation}
\label{ExpABEq}
T^*(a_i)=
a_{i+2}\,
\prod_{k=1}^m
\frac{1+a_{i+3k+2}\,b_{i+3k+1}}{1+a_{i-3k+2}\,b_{i-3k+1}},
\qquad 
T^*(b_i)=b_{i-1}\,
\prod_{k=1}^m
\frac{1+a_{i-3k-2}\,b_{i-3k-1}}{1+a_{i+3k-2}\,b_{i+3k-1}}.
\end{equation}
\end{proposition}

\proof
According to Lemma \ref{TmapDual}, $T=\a \circ \b$.
Combining Lemmas \ref{ExplAlpha} and \ref{ExplBeta},
one obtains the expression:
$$
T^*(a_i)=\frac{\lambda_{i+4}\,a_{i+2}}{\lambda_{i+2}},
\qquad
T^*(b_i)= \frac{\lambda_{i}\,b_{i-1}}{\lambda_{i+2}},
$$
where $\lambda_i$ are as in (\ref{LamB}).
Equation (\ref{RelEq}) then gives
$$
\begin{array}{rclcl}
T^*(x_{i})&=&
\displaystyle
\frac{T^*(a_{i-2})}{T^*(b_{i-2})\,T^*(b_{i-1})}
&=&
\displaystyle
\frac{\lambda_{i+2}\,a_{i}}{\lambda_{i}}\,
\frac{\lambda_{i}}{\lambda_{i-2}\,b_{i-3}}\,
\frac{\lambda_{i+1}}{\lambda_{i-1}\,b_{i-2}}\\[16pt]
&=&
\displaystyle
\frac{a_{i}}{b_{i-2}\,b_{i-3}}\,
\frac{1+b_{i-3}\,a_{i-2}}{1+b_{i-1}\,a_{i}}
&=&
\displaystyle
\frac{a_{i-3}}{b_{i-2}\,b_{i-3}}\,
\frac{1+b_{i-3}\,a_{i-2}}{1+b_{i-1}\,a_{i}}\,
\frac{a_{i}}{a_{i-3}}\\[16pt]
&=&
\displaystyle
x_{i-1}\,
\frac{1-\frac{1}{x_{i-1}\,y_{i-1}}}{1-\frac{1}{x_{i+1}\,y_{i+1}}}\,
\frac{x_i\,y_{i-1}}{x_{i+1}\,y_{i+1}}
&=&
\displaystyle
x_i\,\frac{1-x_{i-1}\,y_{i-1}}{1-x_{i+1}\,y_{i+1}},
\end{array}
$$
and similarly for $y_i$.
We thus proved formula (\ref{ExpXEq}). To prove (\ref{ExpABEq}), one now uses (\ref{xtoa}).
\proofend

\subsection{The Poisson bracket in the $(a,b)$-coordinates}

The explicit formula of the Poisson bracket in the $(a,b)$-coordinates
is more complicated than~(\ref{PB}).
Recall that $n$ is not a multiple of 3 so that we assume $n=3m+1$ or $n=3m+2$.
In both cases the Poisson bracket is given by the same formula.

\begin{proposition}
\label{BrackAB}
The Poisson bracket (\ref{PB}) can be rewritten as follows.
\begin{equation}
\label{PABEq}
\begin{array}{rcl}
\{a_i, a_j\}&=&
\displaystyle
\sum_{k=1}^m
\left(
\delta_{i,j+3k}-\delta_{i,j-3k}
\right)a_i\,a_j,\\[16pt]
\{a_i, b_j\}&=&0,\\[6pt]
\{b_i, b_j\}&=&
\displaystyle
\sum_{k=1}^m
\left(
\delta_{i,j-3k}-\delta_{i,j+3k}
\right)b_i\,b_j.
\end{array}
\end{equation}
\end{proposition}

\proof
One checks using (\ref{atox}) that the brackets between the coordinate
functions $(x_i,y_j)$ coincide with (\ref{PB}).
\proofend

\begin{figure}[hbtp]
\centering
\includegraphics[width=5in]{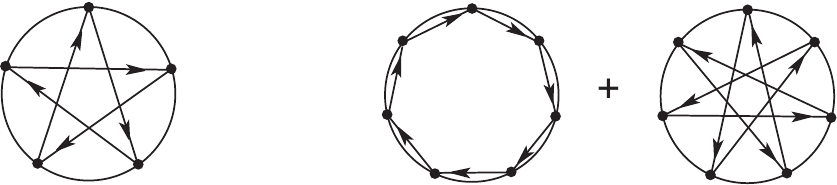}
\caption{The Poisson bracket for $n$=5 and $n=7$.}
\label{Koleso}
\end{figure}


\begin{example}
{\rm
a)
For $n=4$, the bracket is
$$
\{a_i, a_j\}=
\left(
\delta_{i,j+1}-\delta_{i,j-1}
\right)a_i\,a_j
$$
(and with opposite sign for $b$), the other terms vanish.

b)
For $n=5$, the non-zero terms are:
$$
\{a_i, a_j\}=
\left(
\delta_{i,j+2}-\delta_{i,j-2}
\right)a_i\,a_j,
$$
corresponding to the ``pentagram'' in Figure \ref{Koleso}.

c)
For $n=7$, one has:
$$
\{a_i, a_j\}=
\left(
\delta_{i,j+1}-\delta_{i,j-1}
-\delta_{i,j+3}+\delta_{i,j-3}
\right)a_i\,a_j,
$$

d)
For $n=8$, the result is
$$
\{a_i, a_j\}=
\left(
\delta_{i,j+2}-\delta_{i,j-2}
-\delta_{i,j+3}+\delta_{i,j-3}
\right)a_i\,a_j.
$$
}
\end{example}

\section{Monodromy invariants in $(a,b)$-coordinates} \label{monodro}

The $(a,b)$-coordinates are especially well adapted to the computation of the monodromy matrix and the monodromy invariants. Such a computation provides  an alternative deduction of the invariants (\ref{OandE}), independent of \cite{Sch3}.

\subsection{Monodromy matrices} \label{monomat}

Consider the $3\times \infty$ matrix $M$ constructed recurrently as follows: the columns $C_0,C_1,C_2,\dots$ satisfy the relation
\begin{equation} \label{Crec}
C_{i+3}= a_i C_{i+2}+b_i C_{i+1}+C_{i},
\end{equation}
and the initial $3 \times 3$ matrix $(C_0,C_1,C_2)$ is unit. The matrix $M$ contains the monodromy matrices of twisted $n$-gons for all $n$; namely, the following result holds.

\begin{lemma} \label{unimat}
The $3\times 3$ minor $M_n=(C_n,C_{n+1},C_{n+2})$ represents the monodromy matrix of twisted $n$-gons considered as a polynomial function in $a_0,\dots ,a_{n-1},b_0,\dots ,b_{n-1}$.
\end{lemma}

\proof
The recurrence (\ref{Crec}) coincides with (\ref{recur}), see Section \ref{DiffEq}.
It follows that $M_n$ represents the monodromy of twisted $n$-gons in the basis $C_0,C_1,C_2$.
\proofend

Let 
$$
N_j=\left(\begin{array}{ccc}
0&0&1\\
1&0&b_j\\
0&1&a_j
\end{array}\right).
$$
The recurrence (\ref{Crec})  implies the following statement.

\begin{lemma} \label{prod}
One has: $M_n=N_0 N_1\dots N_{n-1}$. In particular,  $\det M_n=1$.
\end{lemma}

To illustrate, the beginning of the matrix $M$ is as follows:
$$
\left(\begin{array}{ccccccc}
1&0&0&1&a_1&a_1a_2+b_2&\dots\\
0&1&0&b_0&b_0a_1+1&b_0a_1a_2+a_2+b_0b_2&\dots\\
0&0&1&a_0&a_0a_1+b_1&a_0a_1a_2+b_1a_2+a_0b_1+1&\dots
\end{array}\right).
$$

The dihedral symmetry $\sigma$, that reverses the orientation of a polygon, replaces the monodromy matrices by their inverses and acts as follows:
$$
\sigma: a_i \mapsto -b_{-i}, \qquad b_i \mapsto -a_{-i};
$$
this follows from rewriting equation (\ref{recur}) as
$$
V_i=-b_i V_{i+1}-a_iV_{i+2}+V_{i+3},
$$
or from Lemma \ref{atox}.\footnote{Since all the sums we are dealing with are cyclic, we slightly abuse the notation and ignore a cyclic shift in the definition of $\sigma$ in the $(a,b)$-coordinates.} 

Consider the rescaling 1-parameter group
$$
\varphi_{\t}: a_i \mapsto e^{\tau} a_i,\qquad b_i \mapsto e^{-\t} b_i.
$$
It follows from Lemma \ref{atox} that the action on the corner invariants is as follows:
$$
x_{i} \mapsto e^{3\t} x_{i},\qquad y_{i} \mapsto e^{-3\t} y_{i}.
$$
Thus our rescaling is essentially the same as the one in (\ref{rescaling}) with $t=e^{3\t}$.

The trace of $M_n$ is a polynomial  $F_n (a_0,\dots, a_{n-1},b_0,\dots,b_{n-1})$. Denote its homogeneous components in $s:=e^{\t}$ by $I_j,\ j=0,\dots ,[n/2]$; these are the monodromy invariants. One has $F_n=\sum I_j$.
The $s$-weight of $I_j$ is given by the formula:
\begin{equation} \label{weights}
w(j) = 3j-k \quad {\rm if}\quad  n=2k,
\qquad  {\rm and}\quad   w(j) = 3j-k+1\quad {\rm if} \quad n=2k+1
\end{equation}
(this will follow from Proposition \ref{abtr} in the next section). For example, $M_4$ is the matrix
$$
\left(\begin{array}{ccc}
a_1&a_1a_2+b_2&a_1a_2a_3+a_3b_2+a_1b_3+1\\
a_1b_0+1&a_1a_2b_0+a_2+b_0b_2&a_1a_2a_3b_0+a_2a_3+a_3b_0b_2+a_1b_0b_3+b_3+b_0\\
a_0a_1+b_1&a_0a_1a_2+a_2b_1+a_0b_2+1&a_0a_1a_2a_3+a_2a_3b_1+a_0a_3b_2+a_0a_1b_3+a_0+a_3+b_1b_3
\end{array}\right)
$$
and 
$$
I_0=b_0b_2+b_1b_3,
\quad 
I_1=a_0+a_1+a_2+a_3+b_0a_1a_2+b_1a_2a_3+b_2a_3a_0+b_3a_0a_1, 
\quad 
I_2=a_0a_1a_2a_3.
$$
Likewise, for $n=5$,
$$
I_0=\sum(b_0+b_0b_2a_3),
\qquad I_1= \sum(a_0a_1+b_0a_1a_2a_3),
\qquad I_2=a_0a_1a_2a_3a_4,
$$
where the sums are cyclic over the indices $0,\dots ,4$.

One also has the second set of monodromy invariants $J_0,\dots ,J_k$ constructed from the inverse monodromy matrix, that is, applying the dihedral involution $\sigma$ to $I_0,\dots ,I_k$.

\subsection{Combinatorics of the monodromy invariants} \label{invcomb}

We now describe the polynomials $I_i,J_i$ and their relation to the monodromy invariants $E_k,O_k$. 

Label the vertices of an oriented regular $n$-gon by $0,1,\dots,n-1$. Consider marking of the vertices by the symbols $a,b$ and $*$ subject to the rule: each marking should be coded by a cyclic word $W$  in symbols $1,2,3$ where $1=a, 2=*\,b, 3=***$. Call such markings {\it admissible}. If $p,q,r$ are the  occurrences of $1,2,3$ in $W$  then $p+2q+3r=n$;
define the weight of $W$ as $p-q$. Given a marking as above, take the product of the respective variables $a_i$ or $b_i$ that occur at vertex $i$; if a vertex is marked by $*$ then it contributes 1 to the product. Denote by $T_j$ the sum of these products  over all markings of  weight $j$. Then $A:=T_k$ is the product of all $a_i$; let $B$ be the product of all $b_i$; here $k=[n/2]$.

\begin{proposition} \label{abtr}
The monodromy invariants $I_j$ coincide with the polynomials $T_j$. One has: 
$$
E_j=\frac{I_{k-j}}{A}\ {\rm for}\  j=1,\dots ,k,\ {\rm and}\  E_n=(-1)^n \frac{B}{A^2}.
$$
$J_j$ are described similarly by the rule $1=b, 2=a\,*, 3=***$, and are similarly related to $O_j$:
$$
O_j=(-1)^{n+j}\frac{J_{k-j}}{B}\ {\rm for}\  j=1,\dots ,k,\ {\rm and}\  O_n=\frac{A}{B^2}.
$$

\end{proposition}

\proof
First, we claim that the trace $F_n$ is invariant under cyclic permutations of the indices $0,1,...,n-1$. 

Indeed, impose the $n$-periodicity condition: $a_{i+n}=a_i, b_{i+n}=b_i$. 
Let $V_i$ be as (\ref{recur}). The matrix $M_n$ takes $(V_0,V_1,V_2)$ to $(V_n,V_{n+1}, V_{n+2})$. Then the matrix $(V_1,V_2,V_3) \to (V_{n+1}, V_{n+2}, V_{n+3})$ is conjugated to $M_n$ and hence has the same trace. This trace is 
$F_n(a_1,b_1,...,a_n,b_n)$, and due to $n$-periodicity, this equals 
$F_n(a_1,b_1,...,a_0,b_0)$. Thus $F_n$ is cyclically invariant. 

Now we argue inductively on $n$. Assume that we know that $I_j=T_j$ for $j=n-2, n-1, n$. Consider $F_{n+1}$.
Given an  admissible labeling of $n-2, n-1$ or $n$-gon, one may insert $***, *\,b$ or $a$ between any two consecutive vertices, respectively, and obtain an admissible labeling of $n+1$-gon. All admissible labeling are thus obtained, possibly, in many different ways.

We claim that $F_{n+1}$  contains the cyclic sums corresponding to all admissible labeling. Indeed, consider an admissible cyclic sum in $F_{n-2}$ corresponding to a labeled $n-2$-gon $L$. This is a cyclic sum of monomials in $a_0,...,b_{n-3}$; these monomials are located in the matrix $M$ on the diagonal  of its minor $M_{n-2}$. By recurrence (\ref{Crec}), the same monomials will appear on the diagonal of $M_{n+1}$, but now they must contribute to a cyclic sum of variables $a_0,...,b_n$. These sums correspond to the labelings of $n+1$-gon that are obtained from $L$ by inserting $***$ between two consecutive vertices.

Likewise, consider a term in $F_{n-1}$, a cyclic sum of monomials in $a_0,...,b_{n-2}$ corresponding to a labeled $n-1$-gon $L$. By (\ref{Crec}), these monomials are to be multiplied by $b_{n-2}, b_{n-1}$ or $b_n$ (depending on whether they appear in the first, second or third row of $M$) and moved 2 units right in the matrix $M$, after which they  contribute to the cyclic sums in $F_{n+1}$. As before, the respective sums correspond to the labelings of $n+1$-gon obtained from $L$ by inserting $*\,b$ between two consecutive vertices. Similarly one deals with a contribution to $F_{n+1}$ from $F_n$: this time, one inserts symbol $a$.

Our next claim is that each admissible term appears in $F_{n+1}$ exactly once.
Suppose not. Using cyclicity, assume this is a monomial $a_n P$ (or, similarly, $b_n P$). Where could a monomial $a_n P$ come from? Only from the the bottom position of the column $C_{n+2}$ (once again, according to recurrence (\ref{Crec})). But then the monomial $P$ appears at least twice in this position, hence in $F_n$, which contradicts our induction assumption. This completes the proof that $I_j=T_j$.

Now let us prove that $E_j=I_{k-j}/A$.  Consider $E_j$ as a function of $x,y$ and switch to the $(a,b)$-coordinates using Lemma \ref{atox}:
$$
x_1=\frac{a_{-1}}{b_{-1}b_0}, \qquad y_1=-\frac{b_0}{a_{-1}a_0}
$$
and its cyclic permutations. Then
$$
y_0x_1y_1=\frac{1}{a_{-2}a_{-1}a_0}
$$
and the cyclic permutations. An admissible monomial in $E_j$ then contributes the factor
$-b_{i}/(a_{i-1}a_i)$ for each singleton $y_{i+1}$ and $1/(a_{i-2}a_{i-1}a_i)$ 
for each triple $y_ix_{i+1}y_{i+1}$. 

Admissibility implies that no index appears twice. Clear denominators by multiplying by $A$, the product of all $a$'s. Then, for each singleton $y_{i+1}$, we get the factor $-b_i$ and empty space $*$ at the previous position $i-1$, because there was $a_{i-1}$ in the denominator and, for each triple $y_ix_{i+1}y_{i+1}$, we get empty spaces $***$ at positions $i-2,i-1,i$. All other, ``free", positions are filled with $a$'s. In other words, the rule $1=a, 2=*\,b, 3=***$ applies. The signs are correct as well, and the result follows. 

Finally, $E_n$ is the product of all $y_{i+1}$, that is, of the terms  $-b_{i}/(a_{i-1}a_i)$. This product equals $(-1)^n B/A^2$.
\proofend

\begin{remark}
{\rm Unlike the invariants $O_k,E_k$, there are no signs involved: all the terms in polynomials $I_i$ are positive.
}
\end{remark}

Similarly to Remark \ref{sense}, the next lemma shows that one can use Proposition \ref{abtr} even if $n$ is a multiple of 3. In particular, this will be useful  in Theorem \ref{closed} in the next section.

\begin{lemma} \label{stilltrue}
If $n$ is a multiple of 3  then the polynomials $I_j, J_j$ of variables $a_0,\dots, b_{n-1}$ are invariant under the action of the group $\R^* \times \R^*$ given in (\ref{action}).
\end{lemma}

\proof Recall that, by Lemma \ref{prod},  $M_n=N_0 N_1\dots N_{n-1}$ where
$$
N_j=\left(\begin{array}{ccc}
0&0&1\\
1&0&b_j\\
0&1&a_j
\end{array}\right).
$$
The action of $\R^* \times \R^*$ on the matrices $N_j$ depends on $j$ mod 3 and is given by the next formulas:
$$
\left(\begin{array}{ccc}
0&0&1\\
1&0&b_0\\
0&1&a_0
\end{array}\right)
\mapsto
\left(\begin{array}{ccc}
0&0&1\\
1&0&\frac{u}{v}\ b_0\\
0&1&u^2v\ a_0
\end{array}\right)
=
\left(\begin{array}{ccc}
1&0&0\\
0&\frac{u}{v}&0\\
0&0&u^2v
\end{array}\right)
\left(\begin{array}{ccc}
0&0&1\\
1&0&b_0\\
0&1&a_0
\end{array}\right)
\left(\begin{array}{ccc}
\frac{v}{u}&0&0\\
0&\frac{1}{u^2v}&0\\
0&0&1
\end{array}\right),
$$
$$
\left(\begin{array}{ccc}
0&0&1\\
1&0&b_1\\
0&1&a_1
\end{array}\right)
\mapsto
\left(\begin{array}{ccc}
0&0&1\\
1&0&uv^2\ b_1\\
0&1&\frac{v}{u}\ a_1
\end{array}\right)
=
\left(\begin{array}{ccc}
1&0&0\\
0&uv^2&0\\
0&0&\frac{v}{u}
\end{array}\right)
\left(\begin{array}{ccc}
0&0&1\\
1&0&b_1\\
0&1&a_1
\end{array}\right)
\left(\begin{array}{ccc}
\frac{1}{uv^2}&0&0\\
0&\frac{u}{v}&0\\
0&0&1
\end{array}\right),
$$
$$
\left(\begin{array}{ccc}
0&0&1\\
1&0&b_2\\
0&1&a_2
\end{array}\right)
\mapsto
\left(\begin{array}{ccc}
0&0&1\\
1&0&\frac{1}{u^2v}\ b_2\\
0&1&\frac{1}{uv^2}\ a_2
\end{array}\right)
=
\left(\begin{array}{ccc}
1&0&0\\
0&\frac{1}{u^2v}&0\\
0&0&\frac{1}{uv^2}
\end{array}\right)
\left(\begin{array}{ccc}
0&0&1\\
1&0&b_2\\
0&1&a_2
\end{array}\right)
\left(\begin{array}{ccc}
u^2v&0&0\\
0&uv^2&0\\
0&0&1
\end{array}\right).
$$
Note that 
$$
\left(\begin{array}{ccc}
\frac{v}{u}&0&0\\
0&\frac{1}{u^2v}&0\\
0&0&1
\end{array}\right)
\left(\begin{array}{ccc}
1&0&0\\
0&uv^2&0\\
0&0&\frac{v}{u}
\end{array}\right)
=\frac{v}{u}\ E,\ 
\left(\begin{array}{ccc}
\frac{1}{uv^2}&0&0\\
0&\frac{u}{v}&0\\
0&0&1
\end{array}\right)
\left(\begin{array}{ccc}
1&0&0\\
0&\frac{1}{u^2v}&0\\
0&0&\frac{1}{uv^2}
\end{array}\right)
=
\frac{1}{uv^2}\ E,
$$
and
$$
\left(\begin{array}{ccc}
u^2v&0&0\\
0&uv^2&0\\
0&0&1
\end{array}\right)
\left(\begin{array}{ccc}
1&0&0\\
0&\frac{u}{v}&0\\
0&0&u^2v
\end{array}\right)
=u^2v\ E,
$$
where $E$ is the unit matrix. Therefore the $\R^* \times \R^*$-action on $M_n$ is as follows:
$$
M_n \mapsto 
\frac{1}{u^2v}
\left(\begin{array}{ccc}
1&0&0\\
0&\frac{u}{v}&0\\
0&0&u^2v
\end{array}\right)
N_0 N_1\dots N_{n-1}
\left(\begin{array}{ccc}
u^2v&0&0\\
0&uv^2&0\\
0&0&1
\end{array}\right)
\sim M_n,
$$
where $\sim$ means ``is conjugated to". It follows that the trace of $M_n$, as a polynomial in $a_0,\dots, b_{n-1}$, is  $\R^* \times \R^*$-invariant, and so are all its homogeneous components.
\proofend

\subsection{Closed polygons} \label{closedpoly}

A closed $n$-gon (as opposed to merely twisted one) is characterized by the condition that $M_n = Id$. This implies that $\sum I_j =3$ (and, of course, $\sum J_j=3$ as well). There are other linear relations on the monodromy invariants which we discovered in computer experiments. All combined, we found five identities.

\begin{theorem} \label{closed}
For a closed $n$-gon, one has:
$$
\sum_{j=0}^k I_j=\sum_{j=0}^k J_j =3,\quad \sum_{j=0}^k w(j) I_j = \sum_{j=0}^k w(j) J_j=0,\quad \sum_{j=0}^k w(j)^2 (I_j-J_j)=0,
$$
where $k=[n/2]$ and $w(j)$ are the weights (\ref{weights}). 
\end{theorem}

\proof
The monodromy $M\in SL(3,\R)$ is a matrix-valued polynomial function of $a_i,b_i$, and  $M(\tau)=M\circ \varphi_{\tau}$ where $\varphi_{\tau}$ is the scaling action
$$
a_i\mapsto e^{\tau} a_i,\quad b_i \mapsto e^{-\tau} b_i.
$$ 
The characterization of ${\cal C}_n$ is $M(0)=Id$. 

Let $e^{\lambda_1}, e^{\lambda_2}, e^{\lambda_3}$ be the eigenvalues of $M(\tau)$ considered as functions of $a_i,b_i,\tau$. Then $\lambda_i=0$ for $\tau=0$ and each $i=1,2,3$, and 
\begin{equation} \label{lambdas}
\lambda_1 + \lambda_2+ \lambda_3=0
\end{equation}
identically. The eigenvalues of $M^{-1}$ are similar with negative $\lambda$s as exponents. The definition of $I(a,b)$ and $J(a,b)$ implies:
\begin{equation} \label{IJ}
e^{\lambda_1} + e^{\lambda_2} + e^{\lambda_3} = \sum e^{\tau w(j)} I_j,\ 
e^{-\lambda_1} + e^{-\lambda_2} + e^{-\lambda_3} = \sum e^{-\tau w(j)} J_j,
\end{equation}
where $w(j)$ are the weights. Setting $\tau=0$ in (\ref{IJ}) we obtain the obvious relations $\sum I_j=\sum J_j =3$. Differentiating (\ref{IJ}) with respect to $\tau$ and setting $\tau=0$, we get
$$
\sum_{i=1}^3 \lambda_i'(0) = \sum w(j) I_j = \sum w(j) J_j.
$$
By (\ref{lambdas}), the left hand side is zero, and we obtain two other relations stated in the theorem.

Differentiate (\ref{IJ}) twice and set $\tau=0$:
$$
\sum_{i=1}^3 \lambda_i''(0) + \lambda_i'(0)^2=\sum w(j)^2 I_j,\ \sum_{i=1}^3 -\lambda_i''(0) + \lambda_i'(0)^2=\sum w(j)^2 J_j.
$$
Subtract and use the fact that $\sum_{i=1}^3 \lambda_i''(0)=0$ (as follows from (\ref{lambdas}) by differentiating in $\tau$ twice) to obtain:
\begin{equation} \label{fifth}
\sum w(j)^2 (I_j-J_j)=0.
\end{equation}
This is the fifth relation of the theorem.
\proofend

\begin{example} \label{closedpenta}
{\rm In the cases $n=4$ and $n=5$, it is easy to solve the equation $M_n=Id$. For $n=4$, the solution is a single point
$$
a_0=a_1=a_2=a_3=1,\qquad
b_0=b_1=b_2=b_3=-1,
$$
and then $I_0=2, I_1=0, I_2=1$. For $n=5$, one has a 2-parameter set of solutions with free parameters $x,y$:
$$
a_0=x,
\quad
a_1=y,
\quad
a_2=-\frac{1+x}{1-xy}, 
\quad
a_3=-(1-xy), 
\quad
a_4=-\frac{1+y}{1-xy},
$$
and $b_i=-a_{i+2}$.  Hence
$$
I_0=J_0=2-z,\ I_1=J_1=1+2z,\ I_2=J_2=-z\ \ {\rm with}\ \ z=\frac{xy(1+x)(1+y)}{1-xy}.
$$
}
\end{example}

\begin{remark}
{\rm ${\cal C\/}_n$ has codimension 8 in ${\cal P\/}_n$, and we have the five  relations of Theorem \ref{closed}. We conjecture that there are no other relations between the monodromy invariants that hold identically on ${\cal C}_n$.
}
\end{remark}

\section{Continuous limit: the Boussinesq equation} \label{Boussinesq}

Since the theory of infinite-dimensional integrable systems
on functional spaces is much more developed than the theory
of discrete integrable systems, it is important to investigate the
$n\to\infty$ ``continuous limit'' of the pentagram map.

It turns out that the continuous limit of $T$ is the classical Boussinesq equation. 
This is quite remarkable since the Boussinesq equation 
and its discrete analogs are thoroughly studied but,
to the best of our knowledge, their geometrical interpretation remained unknown.

We will also show that the Poisson bracket (\ref{PB})
can be viewed as a discrete version of the well-known first Poisson structure
of the Boussinesq equation.

\subsection{Non-degenerate curves and differential operators}\label{Curva}

We understand the continuous limit of a twisted $n$-gon  as
 a smooth parametrized curve $\g:\R\to\RP^2$ with monodromy:
\begin{equation}
 \label{monoGam}
\g(x+1)=M(\g(x)),
\end{equation}
for all $x\in\R$, where $M\in PSL(3,\R)$ is fixed.
The assumption that every three consequtive
points are in general position corresponds to the assumption that
the vectors $\g'(x)$ and $\g''(x)$ are linearly independent for all $x\in\R$.
A curve $\g$ satisfying these conditions is usually called {\it non-degenerate}.

As in the discrete case, we consider classes of projectively equivalent curves.
The continuous analog of the space ${\cal P}_n$, is 
then the space, ${\cal C}$, of parametrized non-degenerate curves in $\RP^2$
up to projective transformations.
The space ${\cal C}$ is very well known in classical projective differential
geometry, see, e.g., \cite{OT} and references therein.

\begin{proposition}
There exists a one-to-one correspondence between ${\cal C}$ and the space of
linear differential operators on $\R$:
\begin{equation}
\label{ContOp}
A=
\left(
\frac{d}{dx}
\right)^3 
+u(x)\,\frac{d}{dx}+v(x),
\end{equation}
where $u$ and $v$ are smooth periodic functions.
\end{proposition}

This statement is classical, but we give here a sketch of the proof.

\proof
A non-degenerate curve $\g(x)$ in $\RP^2$ has a unique lift to $\R^3$,
that we denote by $\G(x)$, satisfying the condition that the determinant of the vectors
$\G(x),\G'(x),\G''(x)$ (the Wronskian) equals 1 for every $x$:
\begin{equation}
 \label{GammaCond}
\left|
\G(x)\,\G'(x)\,\G''(x)
\right|=1.
\end{equation}
The vector $\G'''(x)$ is a linear combination of $\G(x),\G'(x),\G''(x)$
and the condition (\ref{GammaCond}) is equivalent to the fact that
this combination does not depend on $\G''(x)$.
One then obtains:
$$
\G'''(x)
+u(x)\,\G'(x)+v(x)\,\G(x)=0.
$$
Two curves in $\RP^2$ correspond to the same operator if and only if
they are projectively equivalent.

Conversely, every differential operator (\ref{ContOp}) defines a curve in
$\RP^2$.
Indeed, the space of solutions of the differential equation $A\,f=0$ is
3-dimensional.
At any point $x\in\R$, one considers the 2-dimensional subspace of solutions
vanishing at $x$.
This defines a curve in the projectivization of the space dual to the space of solutions.
\proofend

\begin{remark}
{\rm 
It will be convenient to rewrite the above differential operator as a sum of
a skew-symmetric operator and a (zero-order) symmetric operator:
\begin{equation}
 \label{GammaOp}
 A=
\left(
\frac{d}{dx}
\right)^3 
+\frac{1}{2}
\left(
u(x)\,\frac{d}{dx}+\frac{d}{dx}\,u(x)
\right)
+w(x)
\end{equation}
where $w(x)=v(x)-\frac{u'(x)}{2}$.
The pair of functions $(u,w)$ is understood as
the continuous analog of the coordinates $(a_i,b_i)$.
}
\end{remark}

\subsection{Continuous limit of the pentagram map}\label{ConTT}

\begin{figure}[hbtp]
\centering
\includegraphics[width=2in]{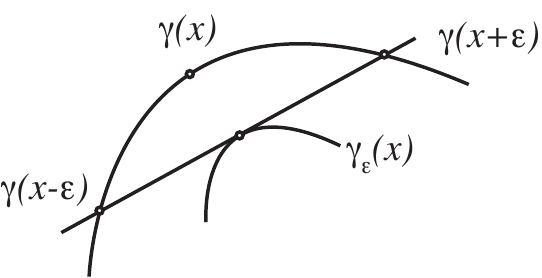}
\caption{Evolution of a non-degenerate curve}
\label{BousFig}
\end{figure}

We are now defining a continuous analog of the map $T$.
The construction is as follows.
Given a non-degenerate curve $\g(x)$, at each point $x$ we draw a
small chord: $\left(\g(x-\eps),\g(x+\eps)\right)$ and obtain a new curve,
$\g_\eps(x)$, as the envelop of these chords, see Figure \ref{BousFig}.

The differential operator (\ref{GammaOp})
corresponding to $\g_\eps(x)$ contains new functions $(u_\eps,w_\eps)$.
We will show that 
$$
u_\eps=u+\eps^2\widetilde{u}+O(\eps^3),
\qquad
w_\eps=w+\eps^2\widetilde{w}+O(\eps^3)
$$
and calculate $\widetilde{u},\widetilde{w}$ explicitly.
We then assume that the functions $u(x)$ and $w(x)$ depend
on an additional parameter $t$ (the ``time'') and define an evolution equation:
$$
\dot{u}=\widetilde{u},
\qquad
\dot{w}=\widetilde{w}
$$
that we undestand as a vector field on the space
of functions (and therefore of operators (\ref{GammaOp})).
Here and below $ \dot{u}$ and $ \dot{w}$ are the partial derivatives in $t$,
the partial derivatives in $x$ will be denoted by $'$.

\begin{theorem}
The continuous limit of the pentagram map $T$ is the following equation:
\begin{equation}
 \label{Bouss}
 \begin{array}{rcl}
  \displaystyle
 \dot{u}&=&w',\\[10pt]
 \displaystyle
\dot{w}&=&
 \displaystyle
 -\frac{u\,u'}{3}-\frac{u'''}{12},
\end{array}
\end{equation}
\end{theorem}

\proof
The (lifted) curve $\G_\eps\subset\R^3$,
corresponding to $\g_\eps(x)$ satisfies the
following conditions:
$$
\begin{array}{rcl}
\left|
\G(x+\eps),\,\G(x-\eps),\,\G_\eps(x)
\right|&=&0,\\[6pt]
\left|
\G_\eps(x),\,\G(x+\eps)-\G(x-\eps),\,\G_\eps'(x)
\right|&=&0.
\end{array}
$$
We assume that the curve $\G_\eps(x)$ is of the form:
$$
\G_\eps=\g+\eps\,A+\eps^2\,B+(\eps^3),
$$
where $A$ and $B$ are some vector-valued functions.
The above conditions easily imply that $A$ is
proportional to $\G$,
while $B$ satisfies:
$$
\begin{array}{rcl}
\left|
\G(x),\,\G'(x),\,B'(x)
\right|&=&0,\\[6pt]
\frac{1}{2}
\left|
\G(x),\,\G'(x),\,\G''(x)
\right|+
\left|
\G(x),\,\G'(x),\,B(x)
\right|&=&0.
\end{array}
$$
It follows that $B=\frac{1}{2}\G''+g\G$, where $g$ is a function.
We proved that the curve
$\G_\eps(x)$ is of the form:
$$
\textstyle
\G_\eps=
\left(1+\eps{}f+\eps^2g\right)\G+
\frac{\eps^2}{2}\,\G''+(\eps^3),
$$
where $f$ and $g$ are some functions.

It remains to find $f$ and $g$ and the corresponding differential operator.
To this end one has to use the normalization condition (\ref{GammaCond}).

\begin{lemma}
The condition (\ref{GammaCond})
implies $f(x)\equiv0$ and $g(x)=\frac{u(x)}{3}$.
\end{lemma}

\proof
A straightforward computation.
\proofend

One has finally the following expression for the lifted curve:
\begin{equation}
\label{GEq}
\G_\eps=
\left(1+\frac{\eps^2}{3}u\right)\G+
\frac{\eps^2}{2}\,\G''+(\eps^3).
\end{equation}
We are ready to find the new functions $u_\eps,v_\eps$ such that
$$
\G_\eps'''(x)
+u_\eps(x)\,\G_\eps'(x)+v_\eps(x)\,\G_\eps(x)=0.
$$
After a straighforward calculation, using the additional
formula
$$
\G^{(V)}=-(2u'+v)\,\G''-(u''+2v'-u^2)\,\G'-(v''-uv)\,\G,
$$
one gets directly from (\ref{GEq}):
$$
u_\eps=u+\eps^2
\left(v'-\frac{u''}{2}\right),
\qquad
v_\eps=v+\eps^2
\left(\frac{v''}{2}-\frac{uu'}{3}-\frac{u'''}{3}\right).
$$
The result follows.
\proofend

\begin{remark}
{\rm
Equation (\ref{Bouss}) is  equivalent to the following differential equation
$$
\ddot{u}= -\frac{\left(u^2\right)''}{6}-\frac{u^{(IV)}}{12},
$$
which is nothing else but the classical Boussinesq equation.
}
\end{remark}

\begin{remark}
{\rm It is not hard to compute that the continuous limit of the scaling symmetry is given by the formula:
$$
u(x)\mapsto u(x),\quad w(x)\mapsto w(x)+t
$$
where $t$ is a constant.
}
\end{remark}

\begin{remark}
{\rm The fact that the continuous limit of the pentagram map is the Boussinesq equation is discovered in \cite{Sch1} (not much details are provided). The computation in \cite{Sch1} is made in an affine chart  $\R^2\subset \RP^2$. In this lift, different from the canonical one (characterized by constant Wronskian), one obtains the curve flow
$$
{\dot \G}=-\frac{1}{3}\frac{[\G',\G''']}{[\G',\G'']}\G' + \frac{1}{2}\G''
$$
where $[.,.]$ is the cross-product; this is equivalent to equation (\ref{GEq}) (we omit a rather tedious verification of this equivalence).
}
\end{remark}

\subsection{The constant Poisson structure}\label{PoiB}

The equation (\ref{Bouss}) is integrable.
In particular, it is Hamiltonian with respect to (two)
Poisson structures on the space of functions $(u,w)$.
We describe here the simplest Poisson structure usually
called the first Poisson structure of the Boussinesq equation.

Consider the space of functionals of the form
$$
H(u,w)=\int_{S^1}h(u,u',\ldots,w,w',\ldots)\,dx,
$$
where $h$ is a polynomial.
The variational derivatives,
$\delta_uH$ and $\delta_wH$, are the smooth functions on $S^1$ given by
the Euler-Lagrange formula, e.g.,
$$
\delta_uH=
\frac{\partial{}h}{\partial{}u}-
\left(
\frac{\partial{}h}{\partial{}u'}
\right)^\prime+
\left(
\frac{\partial{}h}{\partial{}u''}
\right)^{\prime\prime}-+\cdots
$$
and similarly for $\delta_wH$.

\begin{definition}
{\rm
The constant Poisson structure on the space of functionals
is defined by
\begin{equation}
\label{VPBEq}
\{G,H\}=\int_{S^1}
\left(
\delta_uG\left(\delta_wH\right)'
+\delta_wG\left(\delta_uH\right)'
\right)dx.
\end{equation}
Note that the ``functional coordinates'' $\left(u(x),w(x)\right)$ play the role
of Darboux coordinates.
}
\end{definition}

Another way to define the above Poisson structure is as follows.
Given a functional $H$ as above, the Hamiltonian vector field
with Hamiltonian $H$ is given by
\begin{equation}
\label{HamEq}
 \begin{array}{rcl}
  \displaystyle
 \dot{u}&=&
\left(\delta_wH\right)^\prime,\\[6pt]
 \displaystyle
\dot{w}&=&
 \displaystyle
 \left(\delta_uH\right)^\prime.
\end{array}
\end{equation}

The following statement is well known, see, e.g., \cite{FMT}.

\begin{proposition}
\label{BouPropHam}
The function
$$
H=\int_{S^1}\left(
\frac{w^2}{2}-\frac{u^3}{18}-\frac{uu''}{24}
\right)dx
$$
is the Hamiltonian function for the equation (\ref{Bouss}).
\end{proposition}

\proof
Straighforward from (\ref{HamEq}).
\proofend

This statement has many consequences.
For instance, the following functions are the first integrals
of (\ref{Bouss}):
$$
H_1=\int_{S^1}u\,dx,
\qquad
H_2=\int_{S^1}w\,dx,
\qquad
H_3=\int_{S^1}uw\,dx.
$$
Note that the functions $H_1$ and $H_2$ are precisely the
Casimir functions of the structure (\ref{HamEq}).

\subsection{Discretization}\label{DeduS}

Let us now test the inverse procedure of discretization.
Being more ``heuristic'', this procedure will nevertheless
be helpful for understanding the discrete limit of the
classical Poisson structure of the Boussinesq equation.

Given a non-degenerate curve $\g(x)$,
fix an arbitrary point $v_i:=\g(x)$ and, for a small $\eps$,
set $v_{i+1}:=\g(x+\eps)$, etc.
We then have:
$$
v_i=\g(x),
\qquad
v_{i+1}=\g(x+\eps)
\qquad
v_{i+2}=\g(x+2\,\eps)
\qquad
v_{i+3}=\g(x+3\,\eps).
$$
Lifting $\g(x)$ and $(v_i)$ to $\R^3$ so that the difference equation (\ref{recur}) be
satisfied, we are looking for $a$ and $b$ in
$$
\G(x+3\,\eps)=
a(x,\eps)\,\G(x+2\,\eps)+b(x,\eps)\,\G(x+\eps)+\G(x),
$$
as functions of $x$ depending on $\eps$, where $\eps$ is small.

\begin{lemma}
\label{DiscABLem}
Representing $a(x)$ and $b(x)$ as a series in $\eps$:
$$
a(x,\eps)=a_0(x)+\eps\,a_1(x)+\cdots,
\qquad
b(x,\eps)=b_0(x)+\eps\,b_1(x)+\cdots,
$$
one gets for the first four terms:
\begin{equation}
\label{DiscABEq}
\begin{array}{ll}
a_0=3,&b_0=-3,\\[4pt]
a_1=0,&b_1=0,\\[4pt]
a_2=-u(x),&b_2=u(x),\\[6pt]
a_3=-
\frac{7}{4}\,u'(x)-\frac{1}{2}\,w(x),&
b_3=\frac{5}{4}\,u'(x)-\frac{1}{2}\,w(x).
\end{array}
\end{equation}
\end{lemma}

\proof
Straighforward Taylor expansion of the above expression for $\G(x+3\,\eps)$.
\proofend

The Poisson structure (\ref{PB}) can be viwed as a discrete analog of the structure (\ref{VPBEq}) and this is, in fact, the way we guessed it.
Indeed, one has the following two observations.

(1)
Formula (\ref{DiscABEq}) shows that the functions
$\log{a}$ and $\log{b}$ are approximated by linear combinations of $u$ and $w$
and therefore (\ref{PABEq}), a discrete analog of the bracket (\ref{VPBEq}),
should be a constant bracket in the coordinates $(\log{a},\log{b})$.

(2)
Consider the following functionals (linear in the $(a,b)$-coordinates):
$$
F_f(a,b)=\int_{S^1}f(x)\,a(x,\eps)\,dx,
\qquad
G_f(a,b)=\int_{S^1}f(x)\,b(x,\eps)\,dx.
$$
Using (\ref{VPBEq}) and (\ref{DiscABEq}), one obtains:
$$
\{F_f,\,F_g\}=\eps^5\int_{S^1}fg'\,dx+O(\eps^6),
\quad
\{F_f,\,G_g\}=O(\eps^6),
\quad
\{G_f,\,G_g\}=-\eps^5\int_{S^1}fg'\,dx+O(\eps^6).
$$

One immediately derives the following general form of the ``discretized'' Poisson bracket
on the space $\cP_n$:
$$
\{a_i,a_j\}=P_{i,j}\,a_ia_j,
\qquad
\{a_i,b_j\}=0,
\qquad
\{b_i,b_j\}=-P_{i,j}\,b_ib_j,
$$
where $P_{i,j}$ are some constants.
Furthermore, one assumes: $P_{i+k,j+k}=P_{i,j}$ by cyclic invariance.
One then can check that (\ref{PABEq}) is the only Poisson bracket of the
above form preserved by the map $T$.

\section{Discussion}

We finish this paper with  questions and conjectures.

\paragraph*{Second Poisson structure.}

The Poisson structure (\ref{PB}) is a discretization of (\ref{VPBEq})
known as the first Poisson structure of the Boussinesq equation.
We believe that there exists another, second Poisson structure,
compatible with the Poisson structure (\ref{PB}), that allows to
obtain the monodromy invariants (and thus to prove integrability of $T$)
via the standard bi-Hamiltonian procedure.

Note that the Poisson structure usually considered in the
discrete case, cf. \cite{FRS}, is a discrete version of the second
Adler-Gelfand-Dickey bracket. We conjecture that  one can adapt this bracket to the
case of the pentagram map.

\paragraph*{Integrable systems on cluster manifolds.}

The space $\cP_n$ is closely related to cluster manifolds, cf.~\cite{FZ}.
The Poisson bracket (\ref{PB}) is also similar to the 
canonical Poisson bracket on a cluster manifold, cf. \cite{GSV}.

Example \ref{closedpenta} provides a relation of the $(a,b)$-coordinates
to the cluster coordinates.
A twisted pentagon is closed if and only if 
$$
a_0=\frac{a_3+1}{a_1},
\qquad
a_2=\frac{a_1+a_3+1}{a_1\,a_3},
\qquad
a_4=\frac{a_1+1}{a_3},
$$
and $b_i=-a_{i+2}$. Note that $a_{i+5}=a_i$.
This formula coincides with formula of coordinate exchanges for the cluster
manifold of type $A_{2}$, see \cite{FZ}.
The 2-dimensional submanifold of $\cP_5$ with $M=\mathrm{Id}$
is therefore a cluster manifold; the coordinates $(a_1,a_3)$, etc. are the cluster variables.

It would be interesting to compare the coordinate systems on $\cP_n$
naturally arising from our projective geometrical approach with the
canonical cluster coordinates and check whether the Poisson bracket
constructed in this paper coincides with the canonical cluster Poisson bracket.

We think that analogs of the pentagram map may exist for a larger class of
cluster manifolds.

\paragraph*{Polygons inscribed into conics.}
We observed in computer experiments that if a twisted polygon is inscribed into a conic then one has: $E_k=O_k$ for all $k$; the same holds for polygons circumscribed about conics. As of now, this is an open conjecture. Working on this conjecture, we discovered, by computer experiments, a variety of new configuration theorems of projective geometry involving polygons inscribed into conics; these results will be published separately. Let us also mention that twisted $n$-gons inscribed into a conic constitute a coisotropic submanifold of the Poisson manifold ${\cal P}_n$. Dynamical consequences of this observations will be studied elsewhere.

\vskip 1cm

\noindent \textbf{Acknowledgments}.
We are endebted to A. Bobenko, V. Fock, I. Marshall, S. Parmentier, 
M. Semenov-Tian-Shanski and Yu. Suris for stimulating discussions.
V. O. and S. T. are grateful to the Research in Teams program at BIRS. S. T. is also grateful to Max-Planck-Institut in Bonn for its hospitality.
 R. S. and S. T. were partially supported by NSF grants, DMS-0604426  and DMS-0555803, respectively.

\vskip 1cm


Valentin Ovsienko: 
CNRS, Institut Camille Jordan, 
Universit\'e Lyon 1, Villeurbanne Cedex 69622, France, 
ovsienko@math.univ-lyon1.fr

\medskip

Serge Tabachnikov: Pennsylvania State University,
Department of Mathematics University Park, PA 16802, USA, 
tabachni@math.psu.edu

\medskip

Richard Evan Schwartz: 
Department of Mathematics,
Brown University,
Providence, RI 02912, USA, 
res@math.brown.edu

\end{document}